\documentclass{article}

\usepackage{amsfonts}
\usepackage{amsmath}
\usepackage{psfrag} 
\usepackage{epsfig}
\usepackage{subfigure}
\usepackage{graphicx}
\usepackage{graphics}
\usepackage{amssymb}

\usepackage{latexsym}
\usepackage{moreverb}                      
\usepackage{textcomp}

\newenvironment{proof}{\textit{Proof:}\ }{$~\Box$}

\newcommand{\bv}[1]{{\mbox{\boldmath $ #1$}}}
\def\x{{\bv x}}
\def\y{{\bv y}}
\def\p{{\bv p}}
\def\X{{\bv X}}
\def\xsm{{\bv{\scriptstyle x}}}
\def\ysm{{\bv{\scriptstyle y}}}
\def\wi{{\eta}}

\def\Real{\mathbb{R}}  
\def\Znumbers{\mathbb{Z}}  

\newcommand{\grad}[1]{\nabla #1}

\newcommand{\lbeq}[1]{{\label{OR:eq:#1}}}
\newcommand{\be}[1]{\begin{equation} \lbeq{#1}}
\newcommand{\ee}{\end{equation}}
\newcommand{\beno}{\begin{equation*}}
\newcommand{\eeno}{\end{equation*}}
\newcommand{\eq}[1]{{(\ref{OR:eq:#1})}}
\newcommand{\eqtwo}[2]{{(\ref{OR:eq:#1}, \ref{OR:eq:#2})}}

\def\Vector#1{\left(\begin{matrix} #1 \end{matrix}\right)}

\newtheorem{theorem}{Theorem}
\newtheorem{lemma}{Lemma}
\newtheorem{remark}{Remark}

\newcommand{\lbfig}[1]{{\label{OR:fig:#1}}}

\newcommand{\fig}[1]{{Figure~\ref{OR:fig:#1}}}

\begin{document}

\title{Taylor Expansion and Discretization Errors in Gaussian Beam Superposition}
\author{Mohammad Motamed
\thanks{Department of Mathematics, Simon Fraser University, Burnaby, BC V5A 1S6, Canada,
email: {\tt mmotamed@math.sfu.ca}}
\and
Olof Runborg
\thanks{Department of Numerical Analysis, Royal Institute of Technology, 10044 Stockholm, Sweden,
email: {\tt olofr@nada.kth.se}}
}
\date{Feb 4 2010}

\maketitle

\begin{abstract}
The Gaussian beam superposition method is an asymptotic method for computing high frequency wave fields in smoothly varying inhomogeneous media. In this paper we study the accuracy of the Gaussian beam superposition method and derive error estimates related to the discretization of the superposition integral and the Taylor expansion of the phase and amplitude off the center of the beam. We show that in the case of odd order beams, the error is smaller than a simple analysis would indicate because of error cancellation effects between the beams. Since the cancellation happens only when odd order beams are used, there is no remarkable gain in using even order beams. Moreover, applying the
error estimate to the problem with constant speed of propagation, we show that in this case
the local beam width is not a good indicator of accuracy, and there is no direct relation between the error and the beam width. We present numerical examples to verify the error estimates.
\end{abstract}




\section{Introduction}
\label{intro}

Simulation of wave propagation is expensive when the frequency 
of the waves is high. In this case, a large number of grid points are needed to resolve the wave oscillations, and the computational cost to maintain constant accuracy grows algebraically with the frequency. At sufficiently high frequencies, therefore, direct simulations are no longer feasible.

Instead one can use high frequency asymptotic models for wave propagation. The most popular one is geometrical optics, which is obtained when the frequency tends to infinity. The unknowns in geometrical optics are the phase and amplitude which are independent of the frequency and vary on a much coarser scale than the full wave solution. They can therefore be computed at a computational cost independent of the frequency. However, a main drawback of geometrical optics is that the model breaks down at caustics, where geometrical optics rays intersect and the predicted amplitude is unbounded.

Gaussian beams approximation is another high frequency asymptotic model which is valid also at caustics. It was introduced by Popov \cite{Popov}, based on an earlier work of Babic and Pankratova \cite{Babic}. A Gaussian beam is an approximate high frequency solution to the linear wave equation which is concentrated close to a standard ray of geometrical optics, called the central ray of the beam. Although the phase function is real-valued along the central ray, Gaussian beams accept complex-valued phase functions off their central ray. The imaginary part of the phase is chosen such that the solution decays exponentially away from the central ray, maintaining a Gaussian-shaped profile. The main advantage of this construction is that it gives the correct solution also at caustics. It has recently been proved to be beneficial in seismic imaging, \cite{Hill1,Hill2}.

Numerical methods based on Gaussian beams use the superposition principle. Individual beams are computed via ray tracing like equations, where quantities such as the curvature and width of beams are calculated from ordinary differential equations (ODEs) along the central rays, and the contribution of the beams concentrated close to their central rays are determined by Taylor expansion. 
The full wave field is then obtained by a superposition integral over all beams. This integral is replaced by a discrete summation of beams in practical computations.
See for example \cite{Katchalov_Popov,Cerveny_etal,Klimes84,Hill1,Hill2}. Numerical techniques based on both Lagrangian and Eulerian formulations of the problem have been devised \cite{Jin1,Jin2,Qian1,Qian2,Motamed_phd}. For a rigorous mathematical analysis of Gaussian beams we refer to \cite{Ralston}. 

In this paper we derive error estimates for the beam superposition method. We study the discretization error, caused by replacing the superposition integral by the summation of beams, and the error related to Taylor expansion of the phase and amplitude off the center of the beam. Some error estimates for this method have been derived earlier, \cite{Klimesh86,Tanushev}.
We aim to give a more complete picture of the error by also
including the error due to the spreading of the beams, which is related
to the Taylor expansions.
This error is recognized as important in e.g. \cite{Klimesh86}.
It turns out that, in the case of using odd order beams, the error is smaller than a simple analysis would indicate because of error cancellation effects between
the beams. Since the cancellation happens only when odd order beams are used, there is no remarkable gain in using even order beams. Moreover, we show that in the case of constant coefficient equations, i.e. when the speed of propagation is constant, the local beam width is not a good indicator of accuracy, and there is no direct relation between the error and the beams' width. However, this may not be true in the case of varying speed of propagation, where the beam width can be an important factor in the Taylor expansion error.  For other recent results on error estimates see 
\cite{LiuRalston:09,BougachaEtal:09}.

In Section 2, we review the construction of Gaussian beams and the Gaussian beam superposition method. The accuracy of Gaussian beam superposition is studied in Section 3, where the main result is formulated together with numerical examples verifying the obtained error estimates. In Section 4, the proof of the main theorem is given in detail. Finally, in Section 5, we compute the errors analytically in the case of constant coefficient equations and give some remarks on how to select the Gaussian beam parameters.


\section{Gaussian beam superposition method}
\label{gbsum}

Gaussian beams are obtained when the linear wave equation is solved with oscillatory initial or boundary data with an amplitude in the shape of a Gaussian bell. A Gaussian beam is an asymptotic solution concentrated on its central ray in the domain. By the superposition principle for linear equations, such solutions can be added to find the full wave field. The initial/boundary data for beams are obtained such that the wave data at the source is well approximated. In this section, we consider the Helmholtz equation and review the construction of Gaussian beams and their superposition. 

\subsection{Construction of Gaussian beams}
\label{gb_construction}

Consider the Helmholtz equation
$$
\Delta u(\x) + \frac{\omega^2}{c(\x)^2} \, u(\x) = 0, \qquad \x \in \Omega\subset\mathbb{R}^2,
$$
where $\omega \gg 1$ and $c(\x)$ are the frequency and speed of propagation, respectively.
Boundary conditions are given on $\partial\Omega$, which we assume is divided in two parts: one where ingoing waves are specified, and one with outgoing radiation condition, typically at infinity. We call the first, ingoing, part of $\partial\Omega$
the source curve.
We substitute the WKBJ ansatz
\be{WKBJ_ansatz}
u(\x)= e^{i\omega \phi(\xsm)}\sum_{k=0}^\infty A_k(\x)(i\omega)^{-k},
\ee
into the Helmholtz equation. Here, the phase function $\phi$ and the amplitude functions $A_k$ are assumed to be smooth and independent of $\omega$. Equating coefficients of powers of $\omega$ to zero gives us the \emph{eikonal equation} and the \emph{transport equation} for the phase and the first amplitude term in the frequency domain,
$$
 |\grad\phi|=1/c(\x), \qquad 2 \, \grad A_0\cdot\grad\phi + A_0 \, \Delta \phi =0. 
$$
For the remaining amplitude terms, we get
additional transport equations
$$
2\, \grad A_{k+1} \cdot \grad\phi +A_{k+1} \, \Delta \phi + \Delta A_{k}= 0. 
$$
When $\omega$ is large, only the first terms in the WKBJ expansion are significant.
We henceforth denote the high frequency approximation taking $K$ terms in \eq{WKBJ_ansatz} 
by $u_{\rm GO}(\x)$,
\be{Adef}
u_{\rm GO}(\x)=A(\x)e^{i\omega \phi(\xsm)},\qquad
  A(\x) := \sum_{k=0}^{K-1} A_k(\x)(i\omega)^{-k}.
\ee
This approximation is usually called the geometrical optics
approximation, in particular when $K=1$. It introduces an error
of the order $O(\omega^{-K})$.
%

The Gaussian beam approximation has the same form as the geometrical optics
approximation, 
\be{u_gb}
u_{\rm GB}(\x) = A(\x)e^{i\omega \phi(\xsm)},
\ee
where the phase $\phi$ and amplitude terms $A_k$ satisfy the same PDEs. There are,
however, two important differences. First,
while in geometrical optics $\phi$ is globally defined for all rays,
for Gaussian beams it is constructed based on one specific ray (the beam's central ray). Second, in geometrical optics, $\phi$ is real-valued, but in the Gaussian beam construction it is real-valued only on the central ray of the beam. Away from the central ray, it is complex-valued with \emph{positive imaginary part}. The solution will then be exponentially decreasing away from the central ray, maintaining its Gaussian shape. Note that since $\phi$ is complex valued, it actually satisfy
the \emph{complex eikonal equation}, \cite{Magnanini_Talenti1,Magnanini_Talenti2}. Unfortunately, the question of existence and uniqueness of the complex eikonal equation is to a certain extent still open. In particular what precise boundary conditions are well-posed for the above setting is not known.

As in geometrical optics, the Gaussian beam approximation breaks down 
when $\phi(\x)$ becomes non-smooth. 
This is typical for solutions to both the standard and 
the complex eikonal equation.
It happens in general some distance away from the central beam.
On the other hand, away from the beam the solution rapidly goes to zero
and the precise value of the phase is not important. One usually deals with
this problem by multiplying the amplitude with a smooth cut-off function
that is one close to the central ray, and zero for some fixed distance
away from it. In practice, \eq{u_gb} is thus replaced by
$$
u_{\rm GB}(\x) = \varphi(\x)A(\x)e^{i\omega \phi(\xsm)},
$$
where $\varphi(\x)$ is smooth and compactly supported around the central ray.

For a beam starting at point $\x_0$ with direction $\p_0$,
the corresponding central ray satisfies the ray tracing ODEs
\be{rayeqs}
      \frac{d{\x}}{dt} =  c^2(\x)\p, \quad 
    \frac{d{\p}}{dt} = -\frac{\grad c(\x)}{c(\x)},\quad
\x(0)=\x_0,\quad \p(0)=\frac{\p_0}{|\p_0|c(\x_0)},
\ee
with $t$ being the real-valued travel time along the ray. 
If we set $\p=(\cos \theta, \; \sin \theta)^{\top} /c(\x)$ and 
$\x=(x, \; y)^\top$ we can reduce \eq{rayeqs} to
\be{rayeqs2}
\frac{dx}{dt} = c(\x) \cos \theta, \qquad
\frac{dy}{dt} = c(\x) \sin \theta, \qquad
\frac{d\theta}{dt} = c_x(\x) \sin \theta - c_y(\x) \cos \theta.
\ee
The complex-valued $A_k$ and $\phi$ close to the central ray are then approximated by Taylor expansions around the ray,
\begin{align}
\lbeq{Taylor_1}
A_k(\x)&\approx A_k(\x^*) + (\x-\x^*)\cdot\grad A_k(\x^*) + \frac1{2}(\x-\x^*)^{\top} D^2 A_k(\x^*) \, (\x-\x^*) + \cdots,\\
\lbeq{Taylor_2}
\phi(\x)&\approx \phi(\x^*) + (\x-\x^*)\cdot\grad \phi(\x^*)  + \frac1{2}(\x-\x^*)^{\top} D^2 \phi(\x^*) \, (\x-\x^*) + \cdots,
\end{align}
where $\x^*=\x(t)$ for some $t$. The Taylor coefficients
$\phi(\x(t))$, $\grad\phi(\x(t))$, $A_k(\x(t))$, etc. on the central ray can 
be computed. The lowest ones are real on the beam and given directly,
$$
\phi(\x(t))=\phi(\x_0)+t, \qquad \grad\phi(\x(t))=\p(t).
$$
The higher order ones can be obtained by solving ODEs similar to \eq{rayeqs},
and may have a complex part on the beam. 
The most common approximation by far is to take $K=1$ in \eq{Adef}, 
so that $A(\x)=A_0(\x)$, and to
approximate
$A(\x)$ to zeroth order and $\phi(\x)$ to second order. In this case we have, \cite{Klimes84},
\be{APhidef}
A(\x(t))=A(\x_0) \, \left(\frac{c(\x(t))}{c(\x_0)} \, \frac{Q(0)}{Q(t)} \right)^{1/2}, \qquad D^2 \phi(\x(t))=H N H^{\top},
\ee
with
\be{HNdef}
H= \Vector{
\sin \theta   & \cos \theta   \\
-\cos \theta   & \sin \theta  
}, \quad 
N= \Vector{
P/Q   & -c_1/c^2 \\
-c_1/c^2   &  -c_2/c^2 
}, \quad
\Vector{
c_1 \\
c_2 
}=H^{\top} \grad c,
\ee
and the complex-valued scalar functions $P$ and $Q$ satisfy the dynamic ray tracing ODEs
\begin{align}
\label{dynamic1}
\frac{dQ}{dt}&=c^2(\x) \, P, \qquad  Q(0)=Q_0,\\
\label{dynamic2}
\frac{dP}{dt}&=-\frac{c_{xx} \sin^2 \theta -2c_{xy} \sin \theta \cos \theta + c_{yy} \cos^2 \theta}{c(\x)} \, Q, \quad P(0)=P_0. 
\end{align}
The quantities $P$ and $Q$ determine the leading order wavefront curvature and the beam width.
For example, if $\y=\x-\x^*$ is orthogonal to the beam at $\x^*$, then by \eq{u_gb} and \eq{Taylor_2}
$$
  |u_{\rm GB}(\x^*+\y)|\sim\left|e^{i\omega \ysm^\top D^2(\xsm^*)\phi\ysm/2 }\right|
  =e^{-\omega (H^\top\ysm)^\top \Im(N) (H^\top\ysm)/2 }
  =e^{-\omega |\ysm|^2\Im(P/Q)/2},
$$
showing that the effective beam width is proportional to $[\omega\Im(P/Q)/2]^{-1/2}$. 
It can be proved that if $Q_0 \neq 0$ and $\Im (P_0/Q_0)>0$, then $Q(t) \neq 0$ and $\Im (P(t)/Q(t))>0$ along the central ray for all $t>0$, \cite{Popov}. Therefore, by a proper choice of initial data $Q_0$ and $P_0$, each beam will be regular (with finite amplitude at caustics) and concentrate along the central ray. A common choice is $Q_0 > 0$ and $P_0 = i$. 

\subsection{Beam superposition}
\label{beamsum}

Let the source curve be given by $\x_0(s)$ in ${\mathbb R}^2$ parameterized by $s$.
We introduce the
notation $A(\x,s)$, $\phi(\x,s)$ and $\varphi(\x,s)$ for the
amplitude, phase and cut-off of a beam with initial position $\x_0(s)$.
In the Gaussian beam superposition method, the boundary condition on $\x_0(s)$ for the wave field is asymptotically expanded into Gaussian beams, \cite{Klimes84}.
Individual Gaussian beams are computed by solving the ODEs \eq{rayeqs} and (\ref{dynamic1},\ref{dynamic2}). The contributions of the beams concentrated close to their central rays are determined by the approximations \eqtwo{Taylor_1}{Taylor_2} entered in \eq{u_gb}. The wave field is then obtained by the superposition integral over the beams,
\be{superpos_int}
u_s(\x) = \omega^{1/2} \int \varphi(\x,s) \, A(\x,s) \, e^{i\omega\phi(\xsm,s)} ds.
\ee
In practical computations, this integral is replaced by a discrete sum of individual beams, the trapezoidal rule approximation,
\be{wave_field}
u_s^D(\x) = \omega^{1/2} h \sum_{j\in\Znumbers} \varphi(\x,s_j)A(\x,s_j)e^{i\omega\phi(\xsm,s_j)},
\ee
where $h$ is the initial spacing of the beams.
 
The initial conditions for the Taylor coefficient ODEs are chosen such that $u_s^D$ well approximates the exact ingoing boundary data. This can be done in different ways. In particular the initial width of the beams can be varied to give different approximations.
As an example, we consider a plane wave in the $y$-direction as boundary
condition on the $x$-axis, $\x_0(s)=(s,0)$. This will be approximated by a sum of beams starting in the same direction, \cite{Hill1}.
The approximation is based on the relationship
\be{Hill_planewave}
1=\frac{1}{\sqrt{\pi}{\wi}_0}\int e^{-(x-s)^2/{\wi}_0^2} ds=
\sum_j \frac1{\sqrt{\pi}} \, \frac{h}{{\wi}_0} e^{-(x-s_j)^2/{\wi}_0^2} + {\mathcal O}(e^{-({\wi}_0/h)^2}), \quad s_j=jh, 
\ee
with ${\wi}_0$ representing the initial beam widths, see \fig{sum_plane}. 
\begin{figure}[!t]
\centering
\includegraphics[width=0.55\textwidth]{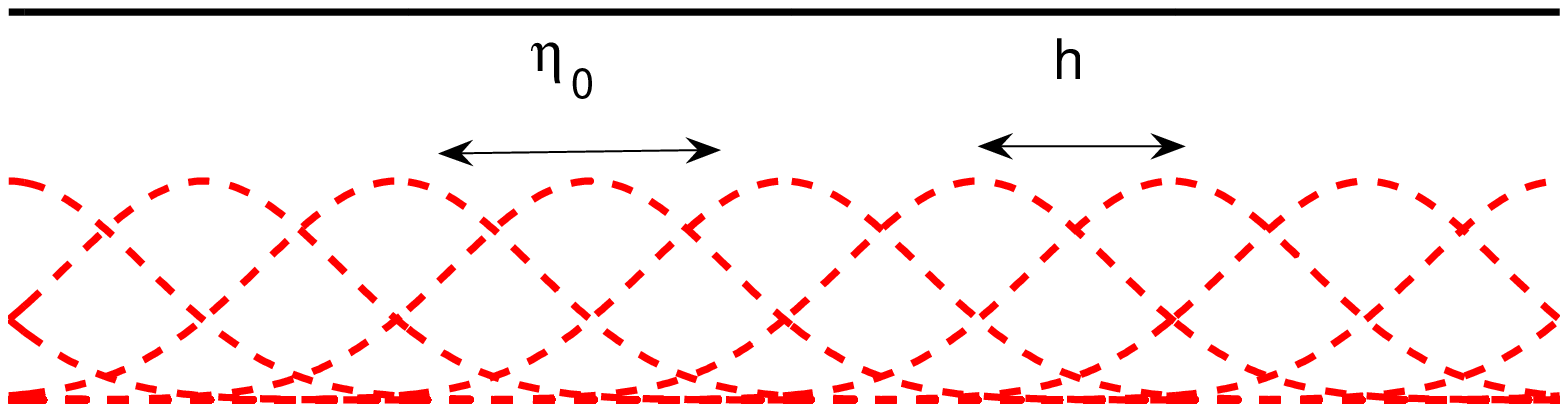}
\begin{picture}(0,0)
\put(-189,20){\line(0,1){23}}
\put(-142,21){\line(0,1){22}}
\put(-100,35){\line(0,1){10}}
\put(-67,35){\line(0,1){10}}
\end{picture}
\caption{The sum of several Gaussian functions is almost constant. A plane wave can therefore be decomposed approximately to a sum of parallel Gaussian beams.}
\lbfig{sum_plane}
\end{figure}
Identifying \eq{Hill_planewave} with \eqtwo{superpos_int}{wave_field}, assuming $\varphi\equiv 1$ we see that
$$
  A(x,0,s)= \frac{1}{\sqrt{\pi\omega} {\wi}_0}, \qquad \phi(x,0,s) = i\frac{(x-s)^2}{\omega {\wi}_0^2}.
$$
To properly choose the initial data, one must take the parameters ${\wi}_0$ and $h$ such that ${{\wi}_0}>h$ by \eq{Hill_planewave}. Then the wave field \eq{wave_field} will produce an accurate plane wave on 
$\x_0(s)=(s,0)$.
The condition ${{\wi}_0}>h$ can be related to the initial data $(P_0,Q_0)$
of the dynamic ray tracing ODEs (\ref{dynamic1},\ref{dynamic2}).
Since beams go in the $y$-direction 
$\theta=\pi/2$ and we have $H=I$ and $\phi_{xx}=P/Q$ by \eq{APhidef} and
\eq{HNdef}. Thus,
$$
  \frac{P_0}{Q_0} = \frac{2i}{\omega {\wi}_0^2}.
$$
Chosing $P_0=i$ we therefore get
$$
h<{\wi}_0= \left( \frac{2 Q_0}{\omega} \right)^{1/2}. 
$$
With this relation between $h$ and $Q_0$ we get an accurate approximation.
In particular we need a spacing $h$ of order $O(1/\sqrt{\omega})$ for a fixed $Q_0$.
Note also that for computational efficiency, $h$ should not be taken much smaller. 
These restrictions were derived for a plane wave but similar
scalings will be necessary also for more general
boundary data.

In what follows, in order to simplify the calculations, we assume that all beams, originating from $\x_0(s)$, shoot out orthogonally.  
We denote by $\X(t,s)$ the location of the center ray
originating in $\x_0(s)$ after time $t$. We further assume that $\phi(\x_0(s),s)=0$.

We make one observation that will be used in the analysis
below. It is well-known that $\X_t\parallel\grad_x\phi$, $\X_t\cdot\grad_x\phi=1$ and
$\X_s\perp\X_t$ under the assumptions made above. Therefore,
since $\phi(\X(t,s),s)=t$, 
\be{phiy0}
   0=\frac{d}{ds}\phi(\X(t,s),s)=\X_s\cdot\grad_x\phi +\phi_s
= \phi_s(\X(t,s),s)
\ee
Hence $\phi_s=0$ everywhere on the central rays.



\section{Accuracy of Gaussian beams summation}
\label{gbsum_accuracy}

In this section we study the accuracy of summation of Gaussian beams. One can distinguish six different types of errors in the approximation:
\begin{enumerate}
\item High frequency approximation.
\item Error in initial data. 
\item Discretization error.
\item Taylor expansion error.
\item Cut-off error.
\item Error in numerical integrators for solving Taylor coefficient ODEs. 
\end{enumerate}
The first error depends on the number of terms used in the WKBJ 
approximation, i.e. the difference $u(\x)-u_{s}(\x)$. For example, for standard beams it is of the order $O(1/\omega)$ since one amplitude term is used, meaning that each beam is a solution to the Helmholtz equation up to order $O(1/\omega)$. The second error represents how well the exact boundary data is approximated by a superposition of Gaussian beams. The third error is caused by replacing the superposition integral by a discrete summation of beams, i.e. $u_s(\x)-u_s^D(\x)$. The fourth error is due to the fact that $A$ and $\phi$ are not computed globally, and only their derivatives on the central beams are computed. One therefore needs to approximate their values around the central beams by Taylor expansions. The fifth error is caused by multiplying the solution by a smooth cut-off function in order to account for possible irregularities away from the central rays. Finally, the last error is the numerical error in solving the ODEs for computing Taylor coefficients. For example the global error in a fourth order Runge-Kutta method is $O(\Delta t ^4)$, with $\Delta t$ being the time-step.

Here, we will only concentrate on the Discretization and Taylor expansion errors. 

\subsection{Motivation and preliminaries}

Let the source be a curve $\x_0(s)$ and assume that we look for the solution along a line $\x=(x,y^*)$. 
We simplify the notation by setting $A_k(\x,s)=A_k(x,s)$, $\phi(\x,s)=\phi(x,s)$,
$\varphi(\x,s)=\varphi(x,s)$ and write
$$
u_s(x) = \omega^{1/2} \int \varphi(x,s) \, A(x,s) \, e^{i\omega\phi(x,s)} ds,
$$
$$
   u_s^D(x) = \omega^{1/2} h \sum_{j\in\Znumbers} \varphi(x,s_j)A(x,s_j)e^{i\omega\phi(x,s_j)}, \qquad s_j=jh.
$$
We now let $X(s)$ denote the location of the center beam on the line $(x,y^*)$ when the initial data is given at $\x_0(s)$. 
Hence, $X(s)$ is implicitly defined by
$$
    \X(t(s),s) = (X(s),y^*),
$$
for some function $t(s)$. \fig{setting} shows the setting for $\x_0(s)=(s,0)$, as an example.
\begin{figure}[!h]
\centering
\psfrag{c}{$\x_0(s)=(s,0)$}
\psfrag{b}{$s$}
\psfrag{a}{$X(s)$}
\psfrag{d}{$\x=(x,y^*)$}
\includegraphics[width=0.6\textwidth]{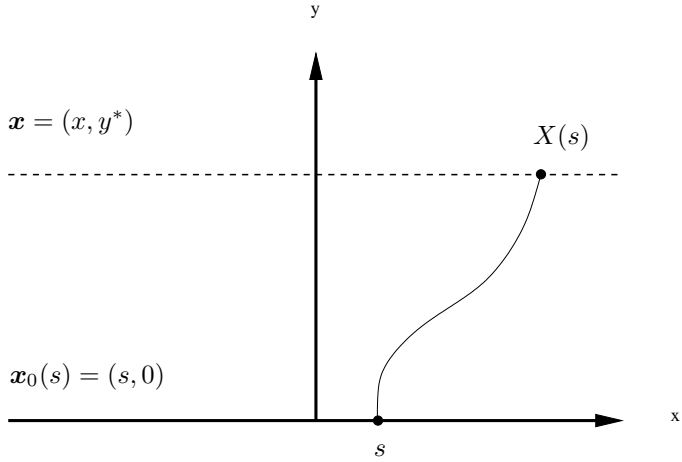}
  \caption{A schematic representation of the initial source and a beam central ray.}
\lbfig{setting}
\end{figure}

We will now explain the approximation of $A(x,s)$ and $\phi(x,s)$
for a $(q+1)$-th order Gaussian beam, complying
with the standard notation that the basic choice $q=0$ is a first order beam.
We then take $K=\lfloor q/2\rfloor+1$ terms in the
WKBJ expansion \eq{Adef} and observe
that the high frequency approximation error will be of the order
\be{hferror}
O(\omega^{-q^*/2})\quad {\rm where}\quad
   q^* = 2(\lfloor q/2\rfloor+1)=\begin{cases}
   q+2, & \text{$q$ even},\\
   q+1, & \text{$q$ odd}.
   \end{cases}
\ee
For the term $A_k(x,s)$ we make a Taylor
expansion up to order $q-2k$ around  $X(s)$,
\be{A_tilde}
   A_k(x,s)\approx \tilde{A}_{k,q-2k}(x,s) := A_k(X(s),s) + \cdots + \frac{(x-X(s))^{q-2k}}{(q-2k)!}\partial_x^{q-2k}A_k(X(s),s),
\ee
with $0\leq k \leq \lfloor q/2\rfloor$. We also set
\be{Aqdef}
\tilde{A}_q(x,s) :=\sum_{k=0}^{\lfloor q/2\rfloor}\tilde{A}_{k,q-2k}(x,s)(i\omega)^{-k}
\approx A(x,s).
\ee
Furthermore, we approximate $\phi(x,s)$ up to level $q+2$,
\be{phi_tilde}
   \phi(x,s)\approx \tilde{\phi}_q(x,s) := \phi(X(s),s) + \cdots + \frac{(x-X(s))^{q+2}}{(q+2)!}\partial_x^{q+2}\phi(X(s),s).
\ee
The approximate Gaussian beam solution is then given by
$$
\tilde{u}_s(x) = \omega^{1/2}
\int \varphi(x,s) \, \tilde{A}_{q}(x,s) \, e^{i\omega \tilde{\phi}_q(x,s)} ds,
$$
$$
   \tilde{u}_s^D(x) = \omega^{1/2} h \sum_{j\in\Znumbers}\varphi(x,s_j) \tilde{A}_{q}(x,s_j)e^{i\omega\tilde{\phi}_q(x,s_j)}.
$$
The reason why the phase is approximated to two orders higher than the amplitude is
to balance the Taylor expansion errors;
the phase error is multiplied by the frequency $\omega$, which
is proportional to one over the beam width squared (cf. \eq{simpleanalysis} below). 
Note that for $q \ge 2$, one needs to take $K>1$ in \eq{Adef}, i.e. to
include more terms in the WKBJ expansion in order to also balance the high frequency approximation error and the Taylor expansion error, 
cf. Remark~\ref{taylremark} below and
the discussion in \cite{Tanushev}. 

Our motivation for considering the Taylor expansion error comes from the following observation. We define the width of the Gaussian beam passing through $(x,y^*)$ as
$$
   \wi(x) := \frac{1}{\sqrt{\omega \, \Im\phi_{xx}(x,X^{-1}(x))}}.
$$
Because of the term $e^{i\omega (x-X(s))^2 \phi_{xx} /2}$ the solution
will be close to zero for $|x-X(s)|>\wi(x)$. A simple error analysis would therefore 
give the following result. Using \eq{Adef} we have
\begin{align*}
  u_s - \tilde{u}_s
  &= (A-\tilde{A}_q)e^{i\omega\tilde{\phi}_q} +  Ae^{i\omega\tilde{\phi}_q}(e^{i\omega(\phi-\tilde{\phi}_q)}-1)\\
  &= \sum_{k=0}^{\lfloor q/2\rfloor}(i\omega)^{-k}(A_k-\tilde{A}_{k,q-2k})e^{i\omega\tilde{\phi}_q} +  Ae^{i\omega\tilde{\phi}_q}
  (e^{i\omega(\phi-\tilde{\phi}_q)}-1)\\
  &=\sum_{k=0}^{\lfloor q/2\rfloor}
  O(\omega^{-k}{\wi}^{q-2k+1})e^{i\omega\tilde{\phi}_q} + 
Ae^{i\omega\tilde{\phi}_q}(e^{iO(\omega {\wi}^{q+3})}-1).
\end{align*}
Hence, since ${\wi}= O(\omega^{-1/2})$ we would have
\be{simpleanalysis}
  u_s - \tilde{u}_s \sim
  \sum_{k=0}^{\lfloor q/2\rfloor}
  O(\omega^{-(q+1)/2}) + 
Ae^{i\omega\tilde{\phi}_q}O(\omega^{-(q+3)/2+1}) \sim   O(\omega^{-(q+1)/2}).
\ee
In particular, for first order beams with $q=0$, the convergence rate in $\omega$ would be half order, i.e. proportional to $1/\sqrt{\omega}$. This is also
what is observed numerically for a single Gaussian beam. However, 
we now consider two numerical examples using superposition of
first order beams to verify the convergence rate for this case.
Since it is difficult to obtain the exact Gaussian beam superposition
solution $u_s$ we here instead compare $\tilde{u}_s$ with first order ($K=1$) 
geometrical optics $u_{\rm GO}$, which is close enough to $u_s$
to verify or refute the sharpness of \eq{simpleanalysis}; 
the high-frequency error in both $u_s$ and first order geometrical optics
is of the order $O(1/\omega)$. Hence, the predicted error of first order
beams, $O(1/\sqrt{\omega})$, would dominate if it was sharp since
$$
   |\tilde{u}_s-u_{\rm GO}| \leq 
   |\tilde{u}_s-u_s| +
   |u_s-u| +
   |u-u_{\rm GO}| \leq C(1/\sqrt{\omega}+1/\omega)\leq C/\sqrt{\omega}.
$$
In the first example, a plane wave generated on the line $y=0$ propagates orthogonally into the computational domain with a variable speed of propagation. \fig{max_Error}a shows the central rays of Gaussian beams, and \fig{max_Error}c shows the absolute value of the Gaussian beams and geometrical optics solutions along the line $y=0.6$, shown in bold in \fig{max_Error}a. \fig{max_Error}e shows the logarithmic scale of the maximum error between the Gaussian beams solution and the geometrical optics solution.
As can be seen, the convergence rate of the error is surprisingly proportional to $\omega^{-1}$, which is half order better than what we expected. 

In the second example, a plane wave generated on the line $x=0$ propagates with an angle of $45^o$ into the computational domain with a variable speed of propagation. The convergence rate of the error, shown in \fig{max_Error}f, is again proportional to 
$\omega^{-1}$. 
\begin{figure}[!h]
  \begin{center}
   \subfigure[]{\includegraphics[width=5.2cm]{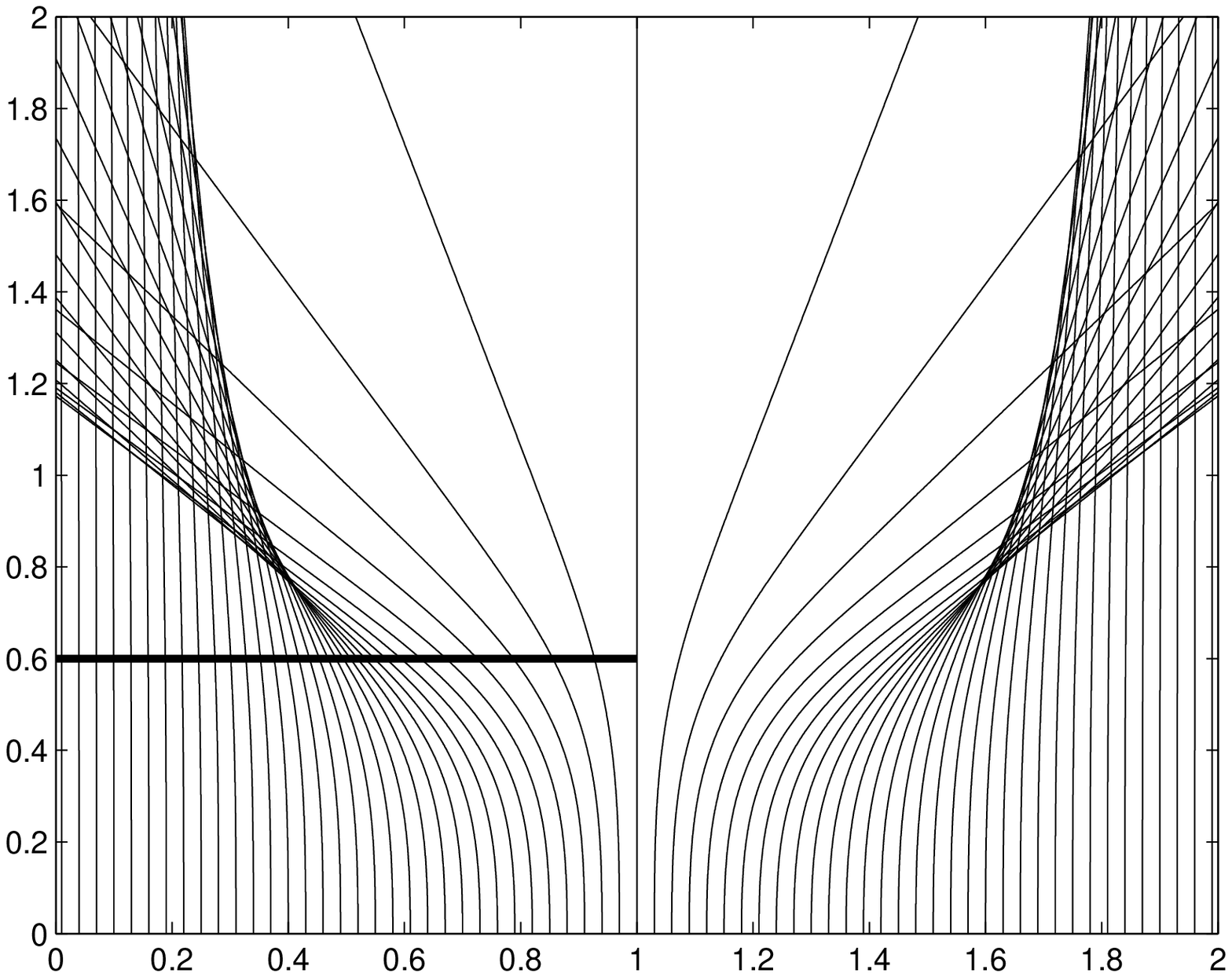}}
   \hskip 9 mm
   \subfigure[]{\includegraphics[width=5.2cm]{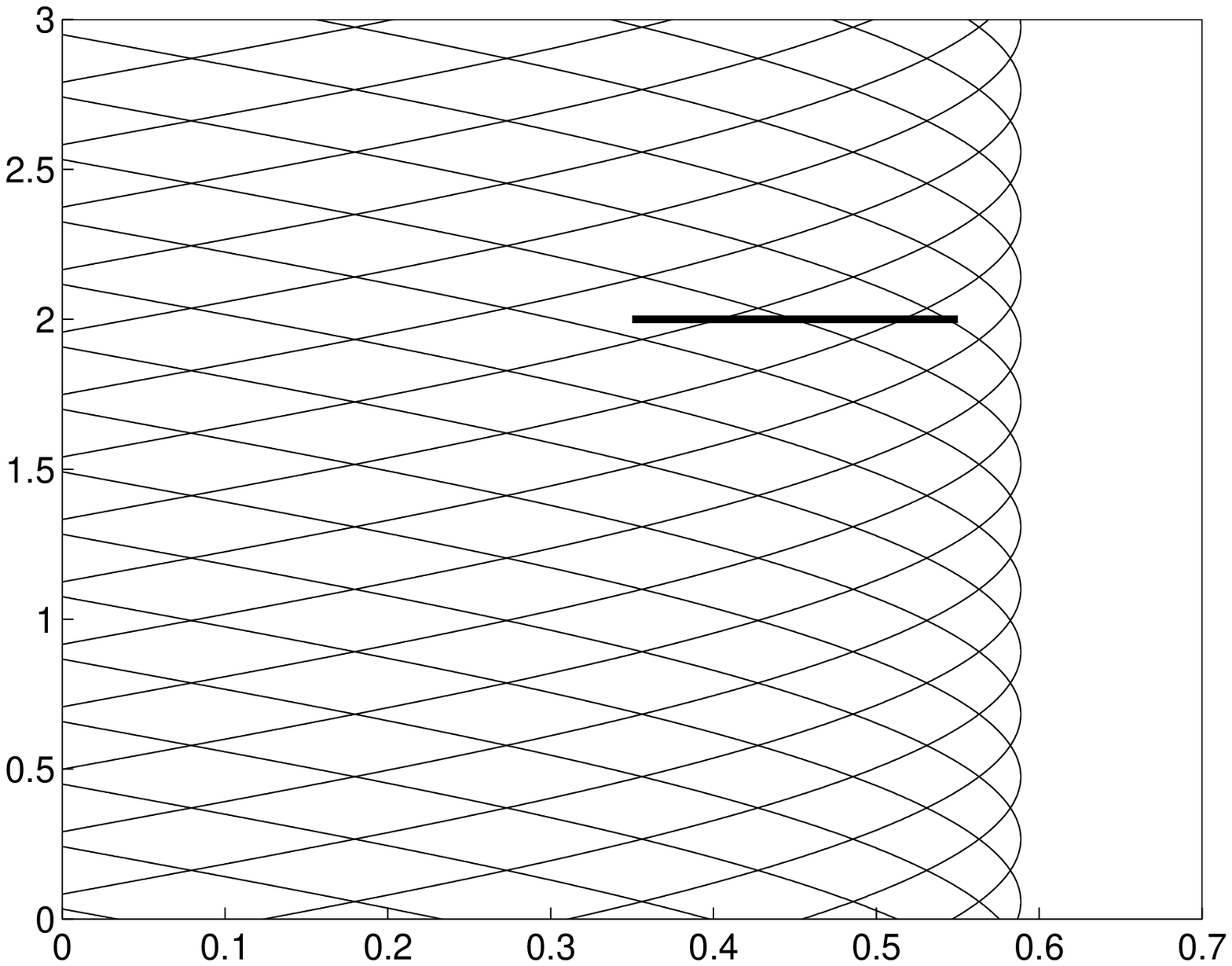}}\\
   \subfigure[]{\includegraphics[width=5.2cm]{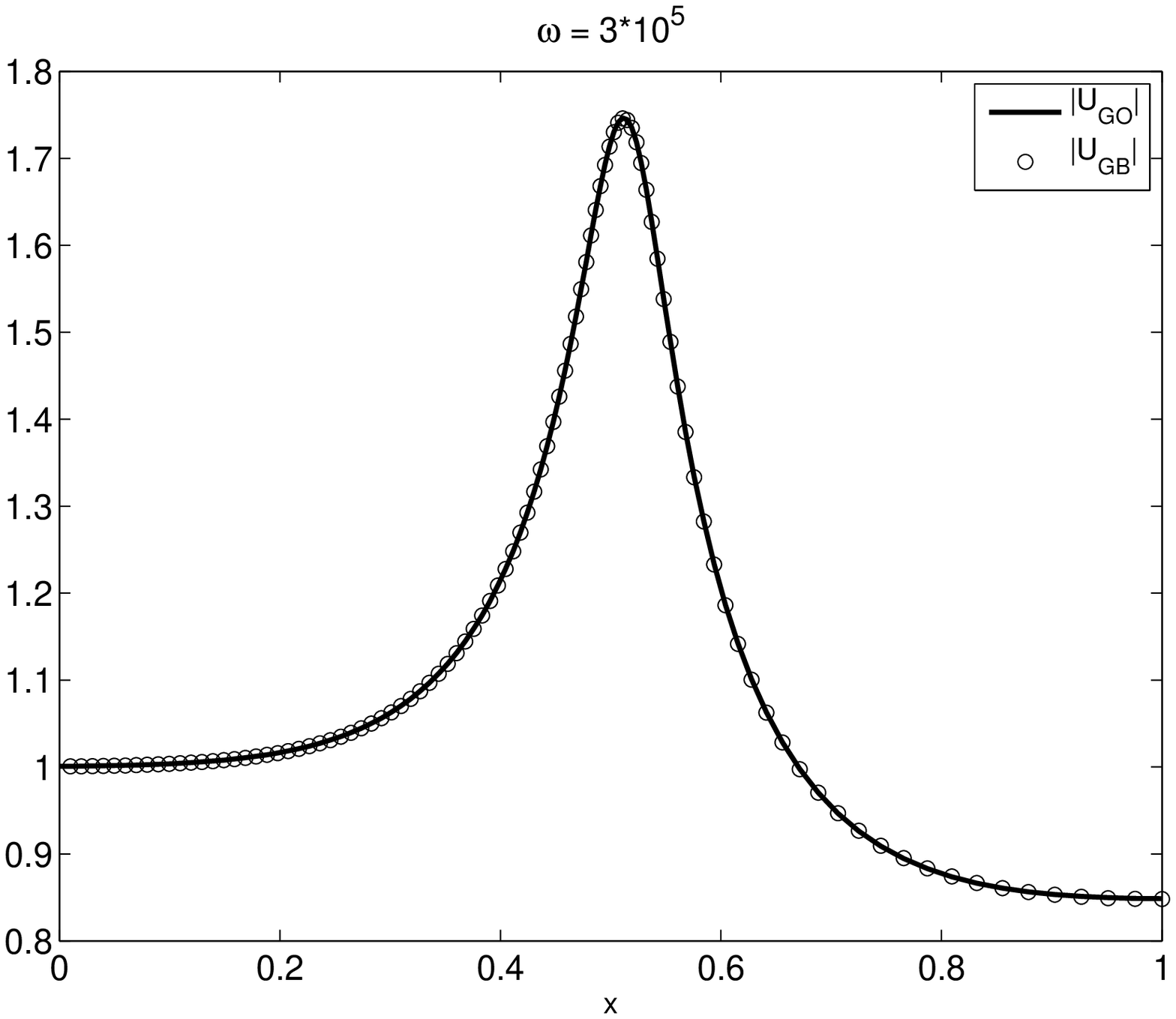}}
   \hskip 9 mm
   \subfigure[]{\includegraphics[width=5.2cm]{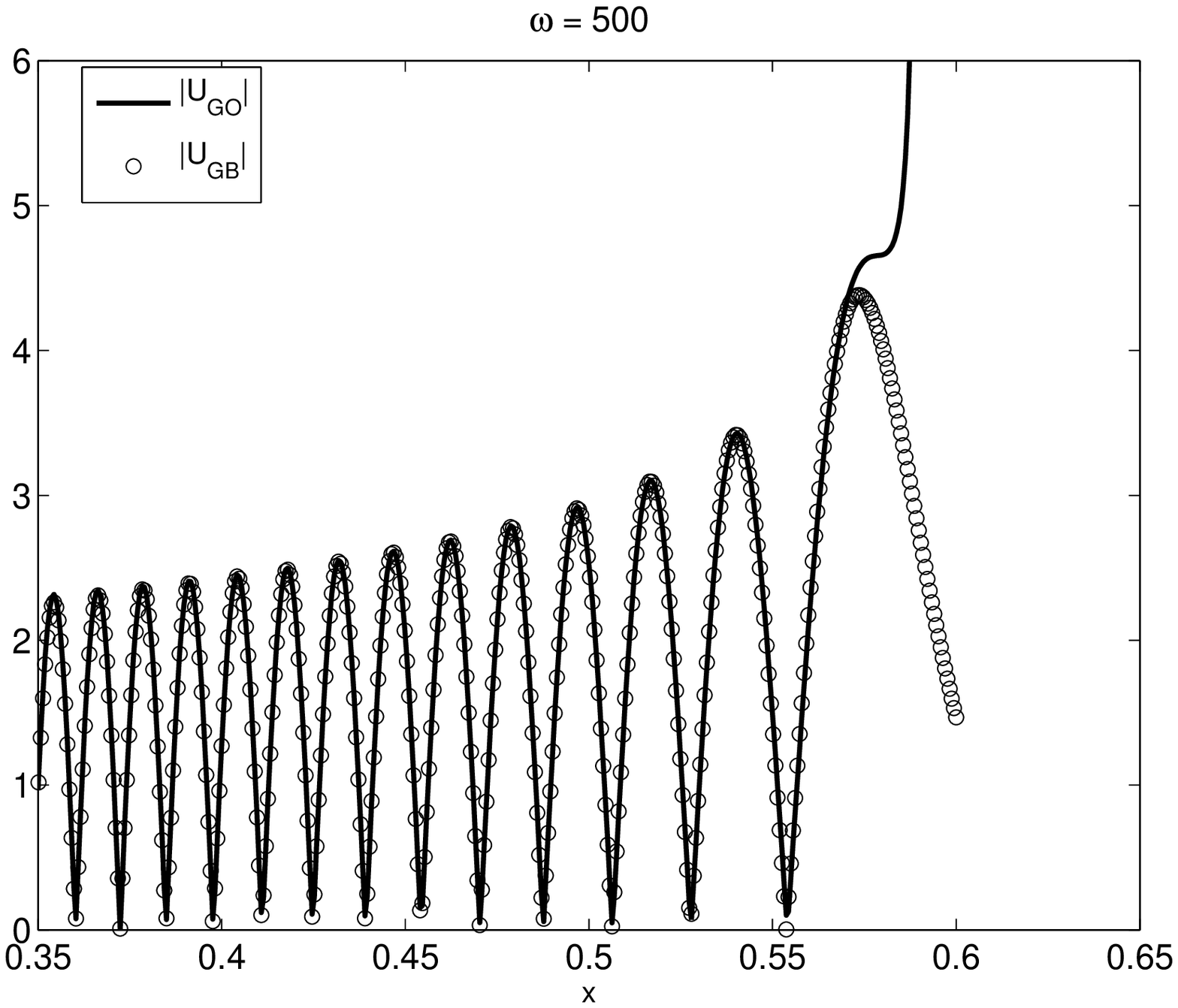}}\\
   \subfigure[]{\includegraphics[width=5.2cm]{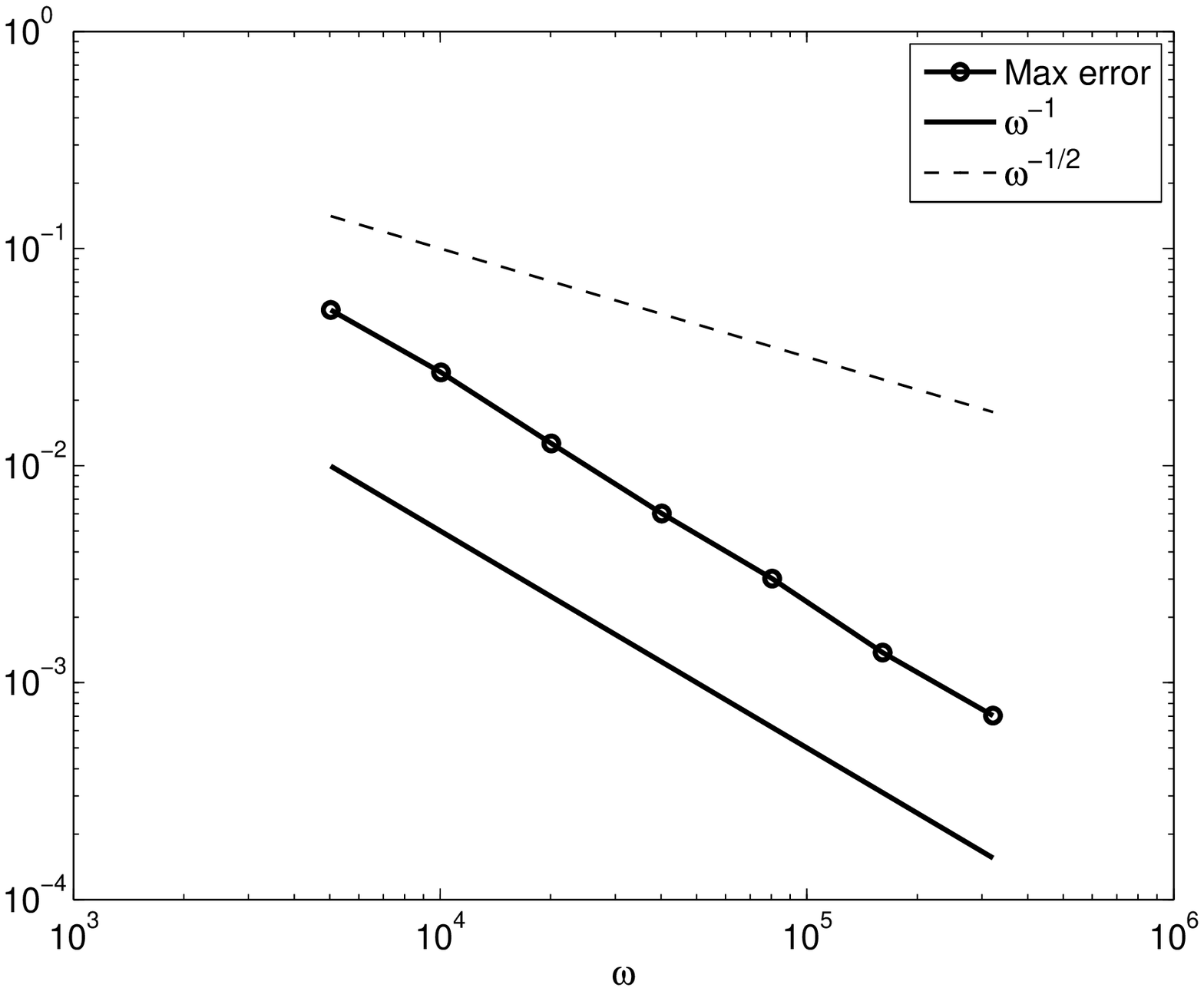}}
   \hskip 9 mm
   \subfigure[]{\includegraphics[width=5.2cm]{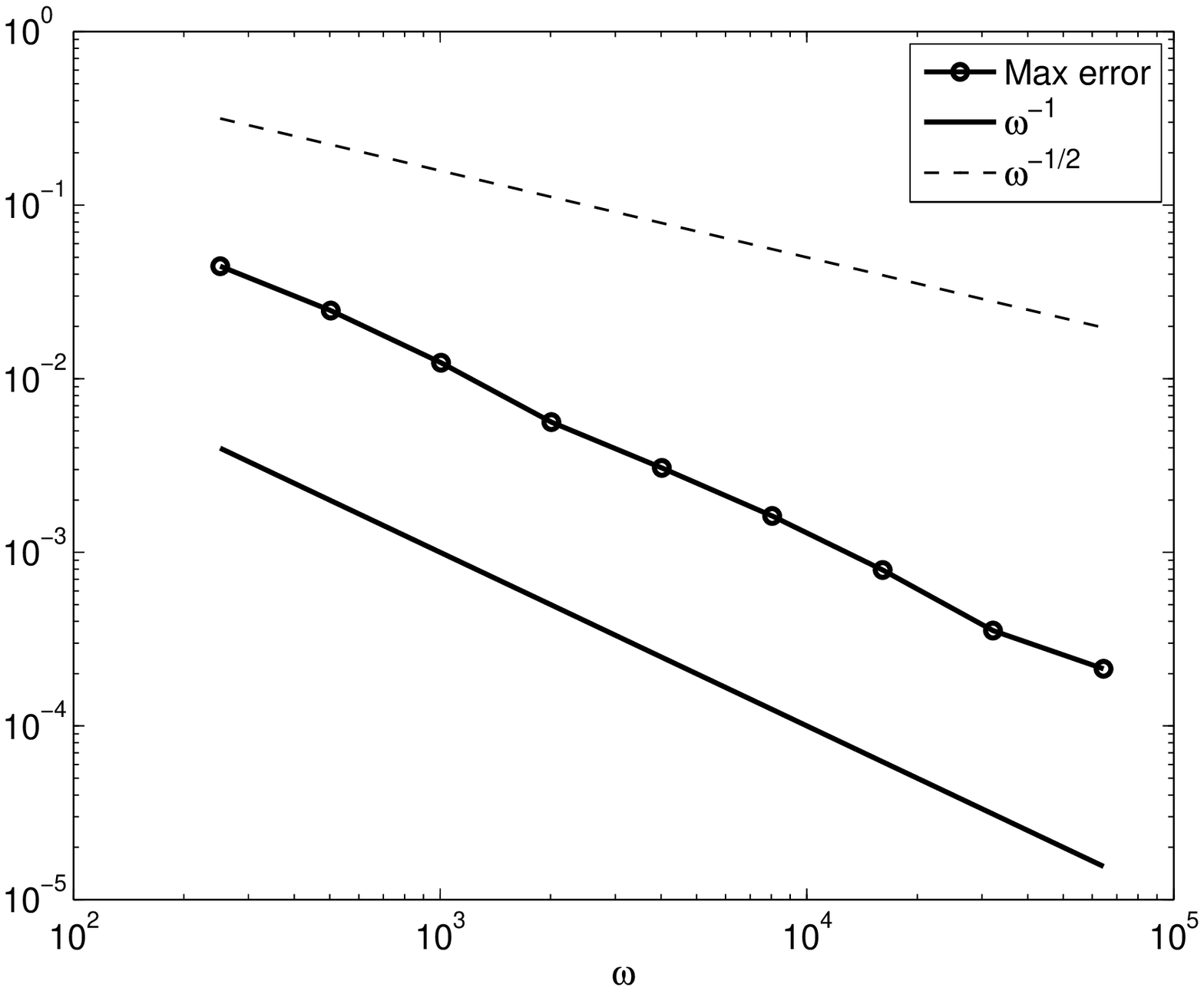}}
  \caption{Left and right top figures show the central rays of Gaussian beams by an initial plane wave on $x$- and $y$-axis, respectively. Middle figures show the absolute value of the Gaussian beams and geometrical optics solutions along the lines $y=0.6$ and $y=2$. Bottom figures show the logarithmic scale of the maximum error between the Gaussian beams solution and the geometrical optics solution. The convergence rate of the maximum error is $\omega^{-1}$.}
\lbfig{max_Error}
  \end{center}
\end{figure}

We will therefore study the Taylor expansion and discretization errors more carefully to describe why it is smaller than what we expected.


\subsection{Main result}

For our results we make the following precise assumptions
\begin{itemize}

\item[({\bf A1})] {\em Smoothness of all coefficients.} We assume $A_k(x,s)\in C_b^{p+q+2-2k}(\Real^2)$, the space of functions with $p+q+2-2k$ continuous and bounded derivatives. Similarly $\phi(x,s)\in C^{p+q+4}$ and $\X(t,s)\in C^{p+1}$, with $p \ge 2$. 

\item[({\bf A2})] {\em Algebraic growth of phase off center beam.} 
For $p_1,p_2\leq p+q+4$, we have for some $\bar{p}$,
$$
\partial_x ^{p_1} \partial_s ^{p_2} \phi(x,s) \leq C (1+|x-X(s)|^{\bar{p}}).
$$
In particular, all derivatives are bounded on the center beam, $x=X(s)$.

\item[({\bf A3})] {\em No caustics.} The derivative $X'(s)$ is bounded away from zero,
$0< c_0\leq X'(s)\leq c_1 <\infty$.

\item[({\bf A4})] {\em Non-degeneracy of each beam.} The imaginary part of $\phi_{xx}$ 
on the beam is
strictly positive and bounded,
\be{posit_imag}
0< c_0 \leq \Im\phi_{xx}(X(s),s)\leq c_1 <\infty.
\ee
Moreover, the frequency is non-vanishing, $\omega>c_2$.
This means that the approximate beams will have a fast decay
off the central beam for high frequencies, and also that the beam width never vanishes or becomes infinite.
The last point is an important feature of Gaussian beams, 
related to the fact that Gaussian beams can approximate the exact field at caustics.

\item[({\bf A5})] {\em Cut-off of fixed size.} We use $\varphi(x,s)=\varphi(X(s)-x)$
with $\varphi\in C^\infty$ such that $\varphi(x)=1$ for $|x|\leq \alpha/2$ 
and $\varphi(x)=0$ for $|x|>\alpha$.
The size of $\alpha$ will be chosen "small enough" depending on $\phi$ but independent of $\omega$. 
\end{itemize}

The error that we want to estimate is given by
$$
E(x) = u_s(x) - \tilde{u}_s^D(x) =  u_s(x) -  \tilde{u}_s(x) + \tilde{u}_s(x) - \tilde{u}_s^D(x) =: E^T + E^D,
$$
where $E^T = u_s(x) -  \tilde{u}_s(x)$ and $E^D = \tilde{u}_s(x) - \tilde{u}_s^D(x)$ represent the Taylor expansion error and the discretization error, respectively.

Then we can show
\begin{theorem}\label{main_thm} 
(Main Theorem) For the $(q+1)$-th order Gaussian beams, we have
$$
|u_s(x) - \tilde{u}^D_s(x)| \leq | E^T| + |E^D|,
$$
where 
\be{th_e1}
| E^T| \leq C \, \omega^{-\frac{q^*}{2}}, \qquad
q^* = \begin{cases}
q+2,& \text{\rm $q$ even},\\
q+1,& \text{\rm $q$ odd},
\end{cases},
\ee
and
\be{th_e2}
|E^D| \leq  C\,\left(\frac{h}{\wi(x)}\right)^p.
\ee
The constants depend on $p$, the initial data, $P_0$ and $Q_0$, for the beams but does not depend on $x$, $\omega$ or $h$. 
For the Taylor expansion error we have
$$
\left| E^T - C^*(x) \, \omega^{1/2}\, {\wi}^{q^* + 1} \right| \leq C' \, \omega^{1/2} \, {\wi}^{q^* + 2},
$$
i.e., the leading order term of the error $E^T$ in $\omega$ is $C^*(x) \, \omega^{1/2} \, {\wi}^{q^* + 1} \sim \omega^{- q^* /2}$, with $C^*(x)$ given by \eq{cstardef},
\eq{bardef}, \eq{e1def} and \eq{e2def}.
\end{theorem}

\begin{remark}
If we take $h<\wi(x)$, the discretization error $E^D$ is typically smaller than the 
Taylor expansion error $E^T$ because of the "spectral" accuracy in \eq{th_e2}.
For first order beams (with $q=0$), the observed
convergence rate is therefore first order in $\omega$, 
which is the same as geometrical optics.
However, $h$ should not be chosen too small for computational complexity reasons. It is also important to note that to balance the error with the error in initial data, $h$ should also relate to the initial beam width ${\wi}_0$.
\end{remark}

\begin{remark}\label{taylremark}
The estimate \eq{th_e1} shows that the Taylor expansion error indeed balances
the high frequency approximation error \eq{hferror}. Moreover, it suggests
that there is no remarkable gain in using even order beams (with an odd $q$);
neither the high frequency nor the Taylor expansion error get a better convergence
rate in $1/\omega$ with these beams.
However, one should note that this is only true in the case of the superposition of beams, where error in adjacent beams cancel. 
If we only have one beam, this does not hold and the
simple error estimate in \eq{simpleanalysis} is sharp. In this case the Taylor
expansion error dominates the high frequency approximation error for even order
beams.
\end{remark}




\section{Proof of main result}


%
%
Before going on to the proof of Theorem \ref{main_thm}, we prove the following utility lemma concerning estimates for the composition of two functions.
\begin{lemma}\label{lemCompo}
Suppose $g_\delta(z)$ belongs to $C^p(\Real)$
for each value of the parameter $\delta$.
If
\be{assump}
   |g^{(k)}_\delta(z)|\leq C_k(1+|z|^q),\qquad 1\leq k \leq p,
\ee
where $C_k$ and $q\geq 0$ are constants independent of $z$ and $\delta$, then there are functions $h_{m,k}\in  C^{p-k}(\Real)$ and constants $C_{m,k}$ independent of $z$ and $\delta$,
such that
\be{Compres}
   \frac{d^k}{dz^k} e^{g_{\delta}(z)} = e^{g_{\delta}(z)} \, \sum_{m=1}^k h_{m,k}(z), 
\qquad \max_{0\leq n\leq p-k} |h^{(n)}_{m,k}(z)|\leq C_{m,k}(1+|z|^{qk}).
\ee
\end{lemma}
\begin{proof}
We show \eq{Compres} by induction. For $k=1$ we have $h_{1,1}=g'_\delta(z)\in C^{p-1}$
and the statement clearly holds. Suppose \eq{Compres} is true for $1 \leq k <p$.
Then
$$
   \frac{d^{k+1}}{dz^{k+1}} e^{g_{\delta}(z)} = e^{g_{\delta}(z)} \, \sum_{m=1}^k h'_{m,k}(z) +g_\delta' (z) \, h_{m,k}(z)=
  e^{g_{\delta}(z)} \, \sum_{m=1}^{k+1} h_{m,k+1}(z)  .
$$
Thus
$$
   h_{m,k+1}(z) = \begin{cases}
h'_{m,k}, & m=1, \\
h'_{m,k}+g_\delta'h_{m-1,k}, & 1<m\leq k, \\
g_\delta'h_{m-1,k}, & m=k+1.
\end{cases}
$$
Using the induction hypothesis, we immediately get that $h_{m,k+1}(z)\in C^{p-k-1}(\Real)$.
Moreover,
$$
  \max_{0\leq n\leq p-k-1}|h^{(n)}_{m,k+1}(z)|
\leq 
\max_{0\leq n\leq p-k-1} |h^{(n+1)}_{m,k}(z)|
+\max_{0\leq n\leq p-k-1}\sum_{j=0}^n c_{jn} |h^{(j)}_{m-1,k}(z)||g^{(n+1-j)}_{\delta}(z)|
$$
The first term is bounded by $C_{1,1}(1+|z|^{qk})$ by assumption, and for
the second term we can estimate
$$
 |h^{(j)}_{m-1,k}(z)||g^{(n+1-j)}_{\delta}(z)|\leq C_{m-1,k}(1+|z|^{qk})\,C_k(1+|z|^q)
\leq C'(1+ |z|^{q(k+1)}),
$$
which proves \eq{Compres}. 
\end{proof}

We can now start with the main proof. We consider each error separately.


\subsection{Taylor expansion error}
\label{ETsec}

The Taylor expansion error is given by
\begin{align}\lbeq{taylerr}
E^T &= u_s(x) - \tilde{u}_s(x) = 
\omega^{1/2} \, \int \varphi(X(s)-x) \, \left( A(x,s) \, e^{i\omega \phi(x,s)}       - \tilde{A}_{q}(x,s) \, e^{i\omega \tilde{\phi}_q(x,s)} \right) ds\\
&=
 \sum_{k=0}^{\lfloor q/2\rfloor} (i\omega)^{-k}\omega^{1/2}\int \varphi(X(s)-x) 
 \left( A_k(x,s) \, e^{i\omega \phi(x,s)}       - \tilde{A}_{k,q-2k}(x,s) \, e^{i\omega \tilde{\phi}_q(x,s)} \right) ds.\nonumber
\end{align}
In this subsection we will start by studying the general
integral approximation
\be{taylred}
\bar{E}^T=\omega^{1/2} \int \varphi(X(s)-x) 
\left( A(x,s) \, e^{i\omega \phi(x,s)}       - \tilde{A}_{q_a}(x,s) \, e^{i\omega \tilde{\phi}_q(x,s)} \right) ds,
\ee
where, with a slight abuse of notation, $A$ and $\tilde{A}_{q_a}$ will represent
one of $A_k$ and $\tilde{A}_{k,q-2k}$ respectively in the sum above. We can
therefore also assume that $q_a\leq q$ and that $q-q_a$ is even.  

Let us denote $X^{-1}(x)$ by $m(x)$ and then, since $X'(s)$ is bounded away from
zero we can use the change of variables
\be{varchange}
  z = \frac{X(s)-x}{\wi(x)} \quad\Rightarrow\quad s = m(x+\wi(x)z).
\ee
We obtain
\begin{align*}
\bar{E}^T & = \omega^{1/2} \, \wi \, \int \varphi({\wi}z) \, \Big( A(x,m(x+{\wi}z)) \, e^{i\omega \phi(x,m(x+{\wi}z))} - 
\\
&
\tilde{A}_{q_a}(x,m(x+{\wi}z)) \, e^{i\omega \tilde{\phi}_q(x,m(x+{\wi}z))} \Big) \, m'(x+{\wi}z) \, dz.
\end{align*}
Now, letting
$$
 D_A(x,s):=A(x,s)-\tilde{A}_{q_a}(x,s), \qquad D_\phi(x,s):=\phi(x,s)-\tilde{\phi}_q(x,s).
$$
and recalling that supp\,$\varphi\subset[-\alpha,\alpha]$,
we can write the integral as
\be{basicint}
\bar{E}^T = \omega^{1/2} \, \wi \, 
\int_{|z|\leq\frac\alpha{w}} \varphi \, \left(D_A + A(e^{i\omega D_\phi}-1) \right) e^{i\omega\tilde{\phi}_q} m' dz.
\ee

We will now approximate the terms in the integral \eq{basicint} by their Taylor expansion. Let us use the shorthand
$$
   \tilde{a}_p(x) = \frac{(-1)^{p}}{p!}\partial_x^{p}A(x,m(x)),\qquad
   \tilde{b}_p(x)  = \left.\frac{1}{p!}\frac{d^p}{dz^p}A(x,m(x+z))\right|_{z=0}.
$$
and
\be{ptildedef}
   \tilde{p}_p(x) = \frac{(-1)^{p}}{p!}\partial_x^{p}\phi(x,m(x)), \qquad
   \tilde{r}_p(x)  = \left.\frac{1}{p!}\frac{d^p}{dz^p}\phi(x,m(x+z))\right|_{z=0}.
\ee
We note that, in this notation
$$
\tilde{A}_{q_a}(x,m(x+z))=\sum_{j=0}^{q_a} \tilde{a}_j(x+z) \, z^j, \qquad \tilde{\phi}_q(x,m(x+z))=\sum_{j=0}^{q+2} \tilde{p}_j(x+z) \, z^j.
$$
Let
$$
   a_1(x) = \tilde{a}_{q_a+1}(x),\qquad
   a_2(x) = \tilde{a}_{q_a+2}(x)+\tilde{a}'_{q_a+1}(x),
$$
$$
   b_1(x) = i\frac{\tilde{p}_{q+3}(x)}{\Im \phi_{xx}(x,m(x))},\qquad
   b_2(x) = i\frac{\tilde{p}_{q+4}(x)+\tilde{p}'_{q+3}(x)}{\Im \phi_{xx}(x,m(x))},
$$
$$
   c_1(x) = \Re\frac{\tilde{r}_{2}(x)}{\Im \phi_{xx}(x,m(x))},\qquad
   c_2(x) = i\frac{\tilde{r}_{3}(x)-\sigma\tilde{p}_{3}(x)}{\Im \phi_{xx}(x,m(x))}.
$$
where $\sigma=1$ for $q=0$ and $\sigma=0$ for $q>0$.
We then approximate
\begin{align*}
   D_A(x,m(x+{\wi}z)) &\approx {\wi}^{q_a+1}\tilde{D}_A(x,z) :=
({\wi}z)^{q_a+1}a_1(x)  + ({\wi}z)^{q_a+2}a_2(x), \\
e^{i\omega D_\phi(x,m(x+{\wi}z))} -1  &\approx {\wi}^{q+1}\tilde{B}(x,z) :=
{\wi}^{q+1}b_1(x)z^{q+3} + {\wi}^{q+2}(b_{2}(x)z^{q+4}+\sigma b_{1}^2(x)z^{2q+6}/2), \\
e^{i\omega\tilde{\phi}_q(x,m(x+{\wi}z))} &\approx \tilde{C}(x,z) =:e^{i\omega\phi(x,m(x))+iz^2c_1(x)-z^2/2}
(1+ c_2(x) \wi z^3).
\end{align*}
The residual terms are denoted
$$
 D_A(x,m(x+{\wi}z))- {\wi}^{q_a+1}\tilde{D}_A(x,z) =: {\wi}^{q_a+3}R_A(x,z),
$$
$$
e^{i\omega D_\phi(x,m(x+{\wi}z))} -1  - {\wi}^{q+1}\tilde{B}(x,z) =:{\wi}^{q+3}R_B(x,z),
$$
$$
e^{i\omega\tilde{\phi}_q(x,m(x+{\wi}z))}- \tilde{C}(x,z)=: {\wi}^2R_C(x,z).
$$
Then we have
\begin{lemma}\label{aux_k0} 
Let the residual terms $R_A$, $R_B$ and $R_C$ be defined as above.
Under assumptions (A1) and (A2), for small enough $\alpha$,
$$
   \left|R_A\right|\leq C|z|^{q_a+3},\qquad
   \left|R_B\right|\leq Ce^{z^2/7},
\qquad
   \left|R_C
\right|\leq Ce^{-z^2/4},
\quad \forall |z|\leq \alpha/\wi,
$$
where  the constant $C$ is independent of $x$, $\omega$ and $z$.
\end{lemma}
\begin{proof}
We note that $\tilde{a}_{q_a}(x+z)$ are the first $q_a$ coefficients in the Taylor expansion
of $A(x+z-x',m(x+z))$ around $x'=0$. Therefore, by Taylor's theorem and assumption (A1)
$$
   \left|D_A(x,m(x+z)) - z^{q_a+1}\tilde{a}_{q_a+1}(x+z)
- z^{q_a+2}\tilde{a}_{q_a+2}(x+z)
\right|\leq C|z|^{q_a+3}.
$$
Expanding the second and third terms around $z=0$ gives the bound for $R_A$.

We now estimate $\omega D_\phi$ in two different ways. By Taylor's theorem as above,
for some $\xi$ with $|\xi-x|\leq {\wi}z$,
$$
|\omega D_\phi(x,m(x+{\wi}z))| \leq \left|\partial^{(q+3)}_x\phi(\xi,m(x+{\wi}z))\right|
\frac{\omega |{\wi}z|^{q+3}}{(q+3)!}\leq
C \omega |{\wi}z|^{q+3}|(1+|{\wi}z|^{\bar{p}}),
$$
where we used the growth condition (A2) for $\phi$ to bound the error term.
Then, for $|z|\leq \alpha/\wi$, and small enough $\alpha$,
\be{pt1}
|\omega D_\phi| \leq C{\wi}^{q+1} |z|^{q+3}, \qquad
|\omega D_\phi| \leq Cz^2\alpha^{q+1}(1+\alpha^{\bar{p}}) \leq \frac{z^2}{8},
\ee
implying
$$
   \left|e^{i\omega D_\phi}-1 - i\omega D_\phi -
   \frac{(i\omega D_\phi)^2}{2}\right|\leq 
   \frac16|\omega D_\phi|^3e^{|\omega D_\phi|}
\leq C{\wi}^{3q+3} |z|^{3q+9} e^{z^2/8}.
$$
Moreover, the same steps as for $D_A$ together with (A2) gives
$$
   \left|D_\phi(x,m(x+z)) - z^{q+3}\tilde{p}_{q+3}(x) -z^{q+4}(\tilde{p}_{q+4}(x)+\tilde{p}'_{q+3}(x))\right|\leq C|z|^{q+5}(1+|z|^{\bar{p}}),
$$
and since $\omega = 1/{\wi}^2\Im\phi_{xx}$, when $|z|\leq\alpha/\wi$,
$$
   \left|i\omega D_\phi(x,m(x+{\wi}z)) - {\wi}^{q+1}z^{q+3}b_1(x) -{\wi}^{q+2}z^{q+4}b_2(x)\right|\leq 
C {\wi}^{q+3}|z|^{q+5}.
$$
Finally, for $q>0$, clearly $|\omega D_\phi|^2\leq C{\wi}^{2q+2}|z|^{2q+6}\leq C {\wi}^{q+3}|z|^{2q+6}$
and for $q=0$ we get 
$$
|(i\omega D_\phi)^2-{\wi}^2b_1^2z^6|=
\frac{|D_\phi^2-({\wi}^3z^3\tilde{p}_3)^2|}{{\wi}^4\Im\phi_{xx}^2}=
\frac{|D_\phi-{\wi}^3z^3\tilde{p}_3||D_\phi+{\wi}^3z^3\tilde{p}_3|}{{\wi}^4\Im\phi_{xx}^2}
\leq C{\wi}^3|z|^7.
$$
Thus,
$$
  |R_B|\leq C {\wi}^{2q} |z|^{3q+9} e^{z^2/8}
+C |z|^{q+5} + (1-\sigma)C |z|^{2q+6} + \sigma C|z|^{2q+7}
\leq C' e^{z^2/7}.
$$

To show the third inequality, we note that since $\phi_s(x,m(x))\equiv 0$ by \eq{phiy0}, we have $\tilde{r}_1(x)=0$. Therefore by Taylor's theorem and assumption (A2), for $q'\geq 2$,
\be{pt2}
  \left|\phi(x,m(x+z)) - \phi(x,m(x))-\sum_{p=2}^{q'}z^p\tilde{r}_p(x)\right|\leq C |z|^{q'+1}(1+|z|^{\bar{p}}).
\ee
Let $v(x,z)= \tilde{\phi}_q(x,m(x+z))-\phi(x,m(x))-z^2\tilde{r}_2(x)$.
Then, by \eq{pt1} and \eq{pt2},
$$
 |v(x,z)| = |\phi(x,m(x+z))-\phi(x,m(x))-D_\phi(x,m(x+z))-z^2\tilde{r}_2(x)|
\leq C|z|^3(1+|z|^{\bar{p}}).
$$
Moreover,
$$
|v(x,z) -z^3(\tilde{r}_3(x) + \sigma \tilde{p}_3(x))|\leq C|z|^4(1+|z|^{\bar{p}}).
$$
As above, if $|z|\leq \alpha/\wi$,
$$
  |e^{i\omega v(x,{\wi}z)} - 1 - {\wi}z^3c_2(x)|\leq |i\omega v(x,{\wi}z) -{\wi}z^3c_2(x)| + \frac12|\omega v|^2e^{|\omega v|}
$$
$$
  \leq C\omega|{\wi}z|^4(1+|{\wi}z|^{\bar{p}}) + \frac12 \left| C\omega|{\wi}z|^3(1+|{\wi}z|^{\bar{p}}) \right|^2e^{C\omega |{\wi}z|^3(1+|{\wi}z|^{\bar{p}})}
$$
$$
  \leq C{\wi}^2|z|^4(1+\alpha^{\bar{p}}) + C{\wi}^2|z|^6(1+\alpha^{\bar{p}})^2e^{Cz^2\alpha(1+\alpha^{\bar{p}})}
\leq C{\wi}^2e^{z^2/4},
$$
for small enough $\alpha$. 
It remains to note that, since $\phi_s(x,m(x))=\Im\phi_x(x,m(x))\equiv 0$,
$$
\Im\tilde{r}_2 = \frac12m'(x)^2\Im\phi_{ss}=\frac12m'(x)\Im\left(\frac{d}{dx}\phi_s - \phi_{sx}\right)
=-\frac12\Im\left(\frac{d}{dx}\phi_x - \phi_{xx}\right)
=\frac12\Im\phi_{xx},
$$
which shows that $i\omega {\wi}^2z^2\tilde{r}_2=iz^2c_1 - z^2/2$.
Therefore,
\begin{align*}
{\wi}^2|R_C| = \left|  e^{i\omega\tilde{\phi}_q(x,m(x+{\wi}z))}- \tilde{C}(x,z) \right| &= \left| \left( e^{i\omega v(x,{\wi}z)} - 1 - {\wi}z^3c_2(x)  \right) \, e^{i \omega \phi(x,m(x)) + i z^2 c_1(x) -z^2/2}    \right| \\
&\leq C {\wi}^2 e^{z^2/4} e^{-z^2/2},
\end{align*}
and the estimate for $|R_C|$ follows. 
\end{proof}

We Taylor expand the remaining quantities in \eq{basicint} and use the assumption (A5) to get
\begin{align*}
\varphi({\wi}z) & \approx 1, \\
A(x,m(x+{\wi}z)) & \approx \tilde{A}(x,z) := A(x,m(x))+{\wi}z\tilde{b}_1(x),\\
m'(x+{\wi}z) & \approx\tilde{m}(x,z) := m'(x) + {\wi}zm''(x).
\end{align*}
It is easy to see that the residual terms for these approximations can
all be bounded by $C{\wi}^2z^2$.
Since these residual terms as well as $R_A$ and $R_B$ above
all grow slower than $\exp(z^2/4)$, we can replace the terms
in the integral in \eq{basicint} by their approximations and control the
error by $O({\wi}^{q_a+3})$,
\be{est1}
   \left| \bar{E}^T- \omega^{1/2} \, \wi \, {\wi}^{q_a+1}\int_{|z|\leq\frac\alpha{\wi}} 
   \tilde{D}_A\tilde{C}\tilde{m}dz
-\omega^{1/2} \, \wi \, {\wi}^{q+1}\int_{|z|\leq\frac\alpha{\wi}} \tilde{A}\tilde{B}
\tilde{C}\tilde{m}dz
\right|\leq C \omega^{1/2} \, \wi \, {\wi}^{q_a+3}.
\ee
Moreover, since the $\tilde{C}(x,z)$ is exponentially small in $\wi$ for 
$|z|\geq\alpha/\wi$
this estimate holds also when taking the integral over all of $\Real$.
We can now compute the leading error terms in $\wi$. The first one, 
$e_1=e_{11}+\wi e_{12}$,
is
\begin{align*}
\lefteqn{\int
   \tilde{D}_A\tilde{C}\tilde{m}dz}\hskip 1 cm 
   \\ &=
   \int (z^{q_a+1}a_1(x)+{\wi}z^{q_a+2}a_2(x))
   e^{i\omega\phi(x,m(x))+iz^2c_1(x)-z^2/2}
(1+ c_2(x) \wi z^3)(m'(x)+{\wi}zm''(x))dz
   \\ &=
   e^{i\omega\phi(x,m(x))}\int z^{q_a+1}a_1m' + {\wi}z^{q_a+2}(
a_2m'+a_1m'' + a_1c_2m'z^2) e^{iz^2c_1(x)-z^2/2}
 dz + O({\wi}^2)   
   \\ &=: e^{i\omega\phi(x,m(x))}(e_{11}+\wi e_{12}) +O({\wi}^2),  
\end{align*}
where
\be{e1def}
  e_{11} = d_{q_a+1}a_1m',
  \qquad
  e_{12} = 
a_2m'd_{q_a+2}+a_1m''d_{q_a+2} + a_1c_2m'd_{q_a+4},
\ee
and
\begin{equation}\label{d_p}
   d_p(x) = \int z^pe^{iz^2c_1(x)-z^2/2}dz = \begin{cases}
N_p (1-2ic_1(x))^{-(p+1)/2},& \text{\rm $p$ even},\\
0,& \text{\rm $p$ odd},
\end{cases}
\end{equation}
with $N_p$ being a constant. We note that $d_p(x)\equiv 0$ when $p$ is odd and it is bounded in $x$ when $p$ is even.
The second term is $e_2=e_{21}+\wi e_{22}$,
\begin{align*}
   \tilde{A}\tilde{B}\tilde{C}\tilde{m}dz 
   =&
   \int 
   [A(x,m(x))+{\wi}z\tilde{b}_1(x)]
   \left[b_1(x)z^{q+3} + \wi\left(b_{2}(x)z^{q+4}+\frac{\sigma}{2} b_{1}^2(x)z^{2q+6}\right)\right]\\
   &\times
   e^{i\omega\phi(x,m(x))+iz^2c_1(x)-z^2/2}
(1+ \wi c_2(x)z^3)(m'(x)+{\wi}zm''(x))dz
   \\ =&\
   e^{i\omega\phi(x,m(x))}
   \int \Bigr[
   Ab_1z^{q+3}m'
  +
   \wi\Bigl( 
\tilde{b}_1b_1z^{q+4}m'
  + A\left(b_{2}z^{q+4}+\frac{\sigma}{2} b_{1}^2z^{2q+6}\right)m'\\
  &
  + Ab_1z^{q+6}c_2m'
  + Ab_1z^{q+4}m''\Bigl)\Bigr]e^{iz^2c_1(x)-z^2/2}dz + O({\wi}^2)
   \\ =:&\ e^{i\omega\phi(x,m(x))}(e_{21}+\wi e_{22}) + O({\wi}^2),
\end{align*}
with
\be{e2def}
  e_{21} = Ab_1d_{q+3}m',
  \qquad
  e_{22} = 
\tilde{b}_1b_1d_{q+4}m'+ A\left(b_{2}d_{q+4}+\frac{\sigma}{2} b_{1}^2d_{2q+6}\right)m'
  + Ab_1d_{q+6}c_2m'
  + Ab_1d_{q+4}m''.
\ee
To find an expression for the leading order error term we now
have to consider four cases depending on $q_a$.
First, if $q_a<q$ then the second term in \eq{est1} is of the same
order or smaller than the right hand side and we can disregard $e_2$.
Second, if $q_a$ is even then $e_{11}=e_{21}=0$ since $d_p=0$ for $p$ odd,
and we gain an additional order in $\wi$. 
Upon also noting that $q-q_a$ is even,
we can therefore write
$$
   \left| \bar{E}^T - \bar{C}(x) \omega^{1/2} \, {\wi}^{\bar{q} + 1} \right|\leq C \omega^{1/2} \, {\wi}^{\bar{q}+2},
$$
where
\be{bardef}
\bar{q}= \begin{cases}
q_a+2,& \text{\rm $q_a$ even},\\
q_a+1,& \text{\rm $q_a$ odd},
\end{cases}
\qquad
\bar{C}(x)=e^{i\omega\phi(x,m(x))}\begin{cases}
e_{11}+e_{21},& \text{\rm $q_a$ odd and $q=q_a$},\\
e_{12}+e_{22},& \text{\rm $q_a$ even and $q=q_a$},\\
e_{11},& \text{\rm $q_a$ odd and $q_a<q$},\\
e_{12},& \text{\rm $q_a$ even and $q_a<q$}.
\end{cases}
\ee
Moreover, $\bar{C}(x)$
is independent of $\omega$ and $h$ and can be bounded by a constant independent of $x$. 

%


%

We now go back to the full Taylor expansion error in \eq{taylerr}
and use the results that were obtained above for \eq{taylred}.
Clearly, all parameters and functions will depend on the term number
$k$, and we indicate this with a subscripted $k$.
Since $q-2k$ is even if and only if $q$ is even, we have
$\bar{q}_k = \bar{q}_0-2k=q^*-2k$ with $q^*$ defined in \eq{th_e1}.
Then,
\begin{align}\lbeq{nonoscerror}
   \left| E^T - \sum_{k=0}^{\lfloor q/2\rfloor}
(i\omega)^{-k}\bar{C}_k(x)\omega^{1/2} \, {\wi}^{\bar{q}_k + 1} \right|
&=
   \left|\sum_{k=0}^{\lfloor q/2\rfloor}(i\omega)^{-k}( \bar{E}^T_k - 
\bar{C}_k(x)\omega^{1/2} \, {\wi}^{\bar{q}_k + 1}) \right|\\
&\leq
  C \omega^{1/2}\sum_{k=0}^{\lfloor q/2\rfloor}\omega^{-k} 
   {\wi}^{\bar{q}_k+2}
=  C \omega^{1/2}\sum_{k=0}^{\lfloor q/2\rfloor}({\wi}^2\omega)^{-k} 
   {\wi}^{q^*+2}.\nonumber
\end{align}
The result therefore follows with
\be{cstardef}
 {C}^*(x) =
  \sum_{k=0}^{\lfloor q/2\rfloor}
(i\omega)^{-k}\bar{C}_k(x)\, {\wi}^{\bar{q}_k - q^*}=
\sum_{k=0}^{\lfloor q/2\rfloor}
(i\omega{\wi}^2)^{-k}\bar{C}_k(x). 
\ee
Since $\omega{\wi}^2=O(1)$,
the leading order term of the error $E^T$ in $\omega$ is $\omega^{-q^* /2}$.

\subsection{Discretization error}

The discretization error is given by
$$
E^D = \tilde{u}_s(x) - \tilde{u}_s^D(x) = \omega^{1/2} \, \int \varphi(X(s)-x) \, \tilde{A}_q(x,s) \, e^{i\omega \tilde{\phi}_q(x,s)} ds - \omega^{1/2} \, h \sum_{j\in\Znumbers} f(j),
$$
with 
$$
f(j) =\varphi(X(s_j)-x) \, \tilde{A}_q(x,s_j) \, e^{i\omega\tilde{\phi}_q(x,s_j)},
$$
for a fixed $x$. The Poisson summation formula gives
$$
\sum_{j\in\Znumbers} f(j) = \sum_{k\in\Znumbers} \hat{f}(k),
$$
where
\begin{align*}
\hat{f}(k) &= \int f(s) e^{-2\pi i s k}ds = \int \varphi(X(sh)-x) \, \tilde{A}_q(x,sh) \, e^{i\omega\tilde{\phi}_q(x,sh)} e^{-2\pi i s k}ds\\
&=\frac1h
  \int \varphi(X(s)-x) \, \tilde{A}_q(x,s) \, e^{i\omega\tilde{\phi}_q(x,s)} e^{-2\pi i s k/h}ds.
\end{align*}
Therefore
$$
E^D = - \omega^{1/2} \, h \sum_{k \ne 0} \hat{f}(k).
$$
Using the change of variables \eq{varchange} we obtain
\be{int_D}
\hat{f}(k) =\frac{\wi}{h}
  \int \varphi({\wi}z) \, \tilde{A}_q(x,m(x+{\wi}z)) \, e^{i\omega\tilde{\phi}_q(x,m(x+{\wi}z))} \, 
e^{-2\pi i m(x+{\wi}z) k/h} \, m'(x+{\wi}z) \, dz.
\ee

We will now show that the integrand functions in \eq{int_D} are smooth, with bounded derivatives. Then the non-stationary phase
lemma can be used to bound $\hat{f}(k)$ since the phase derivative $m'(x)$ never vanishes.

We need
\begin{lemma}\label{aux_k} 
Under assumptions (A1), (A2) and (A4), for $0 \leq \ell \leq p$ and $|z|\leq\alpha / \wi$ with small enough $\alpha$,
\begin{align}
\lbeq{tildeA}
\left|\frac{d^\ell}{dz^\ell} \tilde{A}_q(x,m(x+{\wi}z)) \right| &\leq C,\\ 
\lbeq{e_to_Rephi}
\left|\frac{d^\ell}{dz^\ell} e^{i\omega\tilde{\phi}_q(x,m(x+{\wi}z))} \right| &\leq C' e^{-z^2/5}.
\end{align}
The constants $C$ and $C'$ are independent of of $\ell$, $x$, $\omega $ and $z$.
\end{lemma}
\begin{proof}
For the first inequality we can consider the individual terms in the sum \eq{Aqdef} 
separately. They will each be of the form
$$
\tilde{A}_{k,q-2k}(x,m(x+{\wi}z))=\sum_{j=0}^{q-2k} \tilde{a}_j(x+{\wi}z) \, (\wi z)^j.
$$
where 
$$
   \tilde{a}_j(x) = \frac{(-1)^{j}}{j!}\partial_x^{j}A_{k}(x,m(x)).
$$
By assumption (A1), $\tilde{a}_j \in C_b^{p+2}$ are bounded, for $0 \leq \ell \leq p$, uniformly in $x$ and, after noting that $|{\wi}z|\leq \alpha$ and
that $|\wi|$ is bounded by a constant because of assumption (A4), the 
result \eq{tildeA} follows.

For $\ell=0$ the second inequality is obtained by writing
$$
e^{i\omega\tilde{\phi}_q(x,m(x+{\wi}z))} = \tilde{C}(x,z) + {\wi}^2 R_C, \qquad |   \tilde{C}(x,z)  | = | 1 + {\wi}z^3 c_2(x) | \, e^{-z^2/2}.
$$
Now since $|{\wi}z| \leq \alpha$ and $c_2(x)$ grows algebraically by assumptions (A2), (A3), (A4), and since $\wi$ is bounded by a constant, we have by Lemma~\ref{aux_k0},
\be{case_k0}
|  e^{i\omega\tilde{\phi}_q(x,m(x+{\wi}z))} | \leq  |\tilde{C}(x,z)| + {\wi}^2 | R_C| \leq C e^{-z^2/4}.
\ee

Now consider $1 \leq \ell \leq p$. We write
$$
\tilde{\phi}_q(x,m(x+{\wi}z))=\sum_{j=0}^{q+2} \tilde{p}_j(x+{\wi}z) \, (\wi z)^j,
$$
with $\tilde{p}_j(x)$ defined in \eq{ptildedef}.
Then, since $\tilde{p}_0' + \tilde{p}_1\equiv 0$ by \eq{phiy0}, we have
$$
\frac{d}{dz} \tilde{\phi}_q(x,m(x+{\wi}z))= {\wi}^2 \tilde{p}_1' z + \sum_{j=2}^{q+2} \left( {\wi}^{j+1} z^j \tilde{p}_j'(x+{\wi}z) + j {\wi}^j z^{j-1} \tilde{p}_j(x+{\wi}z) \right),
$$
and therefore
$$
\frac{d}{dz} \left( i \omega \tilde{\phi}_q(x,m(x+{\wi}z)) \right) = \frac{i}{\Im \phi_{xx}(x,m(x))} \, \sum_{j=1}^{q+2} \gamma_j(x+{\wi}z) \, {\wi}^{j-1} z^j,
$$
where
$$
\gamma_j := \tilde{p}_j' + (j+1) \tilde{p}_{j+1}, \quad 1\leq j \leq q+1, \qquad \gamma_{q+2} := \tilde{p}_{q+2}'.
$$
Since the phase derivatives are evaluated on a center beam, $\gamma_j \in C_b^p$ are bounded, for $0 \leq \ell \leq p$, uniformly in $x$ by assumption (A2) and
we therefore have
$$
\left| \frac{d^\ell}{dz^\ell} \left( i \omega \tilde{\phi}_q(x,m(x+{\wi}z)) \right) \right| \leq C_\ell (1+|z|^{q+2}), \qquad 1 \leq \ell \leq p.
$$
Thus, by Lemma~\ref{lemCompo} with $g_{w}=i \omega \tilde{\phi}_q(x,m(x+{\wi}z))$ and $\delta = \wi$, using \eq{Compres} and \eq{case_k0}, the inequality \eq{e_to_Rephi} follows for $1 \leq \ell \leq p$. This completes the proof.
\end{proof} 

The remaining terms in \eq{int_D}, i.e. $\varphi({\wi}z)$ and $m'(x+\wi z)$, are all assumed to be smooth with
derivatives of order up to $p$ bounded uniformly in $x$ by the assumptions (A1) and (A5).
Since the growth in \eq{tildeA} is offset by the rapid decay in \eq{e_to_Rephi},
the above Lemma shows that all $z$-derivatives of the integrand,
$$
   g(x,z):=\varphi \, \tilde{A}_q \, e^{i\omega\tilde{\phi}_q} \, m',
$$
up to order $p$ belongs to $L_1$ and $||\partial^k_zg(x,\cdot)||\leq C_k$ for
$0\leq k \leq p$. The constants
$C_k$ are independent of $x$ and $\omega$.
We can then use the following version of the non-stationary phase lemma.

\begin{lemma}\label{lem3}
Suppose $\psi(z)\in C^{p+1}(\Real)$ with $\psi'(z)\in C_b^{p}(\Real)$ and $\psi'(z)\geq c_0 > 0$.
Moreover, let $\epsilon<\delta<1$ and suppose 
$g(z)\in W^{p,1}$. Then
\be{oscintestimate}
   \left|\int g(z)e^{-i\psi(\delta z)/\varepsilon}dz\right|
\leq 
C ||g||_{W^{p,1}}
\left(\frac{\varepsilon}{\delta}\right)^p,
\ee
where $C$ depends on $\psi(x)$ and $p$, but not on $g(z)$, $\delta$ and $\varepsilon$.
\end{lemma}
For the proof we refer to \cite{Hormander}. It is an easy adaptation of the proof of theorem 7.7.1.

Taking $\psi$ as $2\pi m(x+\cdot)$, $\delta$ as $\wi$ and
$\varepsilon$ as $h/k$ we can apply this to \eq{int_D},
$$
  |\hat{f}(k)|=\frac{\wi}{h}\left|\int g(x,z)e^{-2\pi im(x+{\wi}z)k/h}dz\right|
\leq C \frac{\wi}{h} ||g(x,\cdot)||_{W^{p,1}}
\left(\frac{h}{k \wi}\right)^p.
$$
Consequently,
$$
   \left|\sum_{k\neq 0}\hat{f}(k)\right|\leq
C \frac{\wi}{h} ||g(x,\cdot)||_{W^{p,1}}\sum_{k\neq 0}
\left(\frac{h}{k \wi}\right)^p
\leq C \frac{\wi}{h} 
\left(\frac{h}{\wi}\right)^p
\sum_{k=1}^\infty k^{-p}
\leq C' \frac{\wi}{h} 
\left(\frac{h}{\wi}\right)^p.
$$
Thus since by the assumptions (A1) $p \geq 2$,
$$
|E^D| =  \omega^{1/2} \, h \, |\sum_{k\neq 0} \hat{f}(k)| \leq C'  \omega^{1/2} \, \wi \left(\frac{h}{\wi}\right)^p.
$$
Together with \eq{nonoscerror} this shows the theorem.



\section{Constant coefficient equations}
\label{const_coeff}

It is often claimed that the beam width is important in the accuracy of Gaussian beams, because for wide beams the Taylor expansion error should be large. See for example \cite{Cerveny_etal,Hill1}. We therefore in this section consider the constant coefficient Helmholtz equation, with the speed of propagation $c(\x) \equiv 1$, for which exact Gaussian beam solutions and the Taylor expansion error $|E^T|$ can be computed. We investigate the importance of the beam width on Taylor error in this particular case. Our conclusion is that the local beam width is not a good indicator of accuracy, and there is no direct relation between the error and the beams' width.
We show the main steps of the derivation and the final expression
for $C^*(x)$ and the leading relative error terms below.
For more details we refer to \cite{Motamed_phd}.

We consider first order beams where $q=0$. These only contain one
amplitude term $A_0(x)$ which we for simplicity call $A(x)$ here.
The source curve will be denoted
$\x_0(s)=(s,y_0(s))$ and we assume all beams originating from $\x_0$ shoot out orthogonally. Therefore $\theta_0(s)=\frac{\pi}{2} + \tan^{-1}(y_0'(s))$. In the constant coefficient case the central ray $\Omega$ 
is a straight line.
With $x(0)=x_0(s)=s$, $y(0)=y_0(s)$ and $\theta(0)=\theta_0(s)$, we get from \eq{rayeqs2} at $y=y^*$, 
\begin{gather}
\label{tet_const}
\theta (t(s)) = \frac{\pi}{2} + \tan^{-1}(y_0'(s)), \\ 
\label{X_const}
x(t(s))=X(s)=s-y_0'(s)\, (y^*-y_0(s)), \\
\label{t_const}
 t(s)=\left( (X(s)-s)^2 + (y^*-y_0(s))^2\right)^{1/2}.
\end{gather}
Here we will only compute the error at $\x=(0,y^*)$. For this point, let $s^*:=m(0)=X^{-1}(0)$. To simplify the calculations, and without loss of generality, we assume $y_0(s^*)=y_0'(s^*)=0$. Therefore, by (\ref{tet_const})-(\ref{t_const}), the central ray starting at $\x_0(s^*)$ will lie on the $y$-axis, and we have $s^*=X(s^*)=0$ and $t(s^*)=y^*$. See \fig{circ}.
\begin{figure}[!h]
\begin{center}
\psfrag{a}{$y=y^*$}
\psfrag{b}{$(X(s^*),y^*)$}
\psfrag{x}{$x$}
\psfrag{y}{$y$}
\psfrag{s}{$s^*$}
\includegraphics[width=0.5\textwidth]{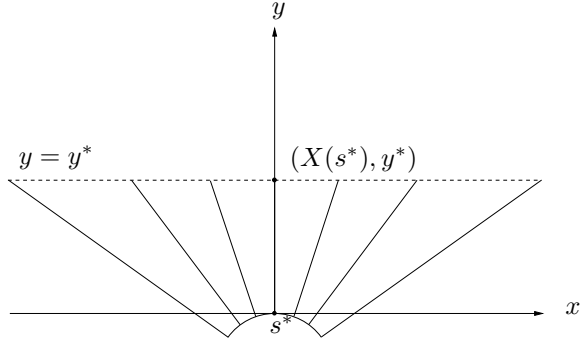}
\caption{A schematic representation of the initial source and central beam rays which are straight lines.}
\lbfig{circ}
\end{center}
\end{figure}

Assuming the initial phase on $\x_0(s)$ to be zero, $\phi(\x_0)=0$, we also get
\be{phi-phix}
\phi(X(s),s)=t(s), \qquad \phi(0,m(0))=y^*.
\ee

To obtain ODEs for higher order Taylor coefficients, we introduce the orthogonal ray-centered coordinates $t,n$, where $n$ is the axis perpendicular to the ray at point $t$ with the origin on the ray.  In this coordinate system, $\phi(t,n=0)$ and $A(t,n=0)$ correspond to $\phi(X(s),s)$ and $A(X(s),s)$ in the Cartesian coordinate, respectively. The eikonal equation and transport equation in the ray-centered coordinates read
\be{eikonal_rayc}
\phi_t^2 + \phi_n^2 = 1,
\ee
\be{transport_rayc}
2\grad A\cdot\grad\phi + A\Delta \phi =0, \qquad \grad \phi = (\phi_t \, \, \phi_n)^{\top}.
\ee
Set $\phi ^{(j)}(t):=\partial_n^j \phi (t,n=0)$ and $A ^{(j)}(t):=\partial_n^j A (t,n=0)$, with $j=0,1,2,\dotsc$. We first note that by \eq{phi-phix}, 
$$
\phi^{(0)}(t)=t, \qquad \partial_t \phi (t,n=0)=1, \qquad \partial_t^j \phi (t,n=0)=0, \quad j=2,3,\dotsc.
$$
Moreover, by \eq{eikonal_rayc} and \eq{transport_rayc} and taking several of their derivatives with respect to $t$ and $n$, 
\begin{gather*}
\phi^{(1)}(t)=0, \quad \partial_t \partial_n \phi(t,n=0)=0, \quad \partial_t \partial_n^2 \phi(t,n=0)=-{\phi^{(2)}}^2(t),\\
\partial_t \partial_n^3 \phi(t,n=0)=0, \quad \partial_t^2 \partial_n \phi(t,n=0)=0, \quad \partial_t^3 \partial_n \phi(t,n=0)=0,\\
\partial_t^2 \partial_n^2 \phi(t,n=0)=2 {\phi^{(2)}}^3(t), \quad \partial_t A(t,n=0)=-\frac1{2} A^{(0)}(t)\, \phi^{(2)}(t),\\
\partial_t^2 A(t,n=0)=\frac3{4} A^{(0)}(t)\, {\phi^{(2)}}^2(t), \quad \partial_t \partial_n A(t,n=0)=0.
\end{gather*}
Now, let
\be{phi-approx}
\phi(t,n)\approx t+\frac{n^2}{2} \phi^{(2)}(t)+\frac{n^3}{6} \phi^{(3)}(t)+\frac{n^4}{24} \phi^{(4)}(t),
\ee
and
\be{A-approx}
A(t,n)\approx A^{(0)}(t)+n A^{(1)}(t)+\frac{n^2}{2} A^{(2)}(t).
\ee
Putting \eq{phi-approx} and \eq{A-approx} into \eq{eikonal_rayc} and \eq{transport_rayc}, we obtain the following ODEs for the Taylor coefficients,
\begin{gather*}
\label{phi2}
\frac{d}{dt} \phi^{(2)} + {\phi^{(2)}}^2 = 0, \\
\label{phi3}
\frac{d}{dt} \phi^{(3)} + 3  \phi^{(2)} \phi^{(3)} =0,\\
\label{phi4}
\frac{d}{dt} \phi^{(4)} + 4 \phi^{(2)} \phi^{(4)} + 3 {\phi^{(2)}}^4 + 3 {\phi^{(3)}}^2=0,\\
\label{A0}
\frac{d}{dt} A^{(0)} + \frac1{2}\phi^{(2)}  A^{(0)}= 0,\\
\label{A1}
\frac{d}{dt} A^{(1)} + \frac3{2}\phi^{(2)}  A^{(1)} + \frac1{2}\phi^{(3)}  A^{(0)} = 0,\\
\label{A2}
\frac{d}{dt} A^{(2)} + \frac5{2}\phi^{(2)}  A^{(2)} + 2\phi^{(3)}  A^{(1)} + \frac1{2}\phi^{(4)}  A^{(0)} + \frac3{2} {\phi^{(2)}}^3  A^{(0)}= 0
\end{gather*}
We then solve these ODEs with $A^{(0)}(0) = 1$ and zero initial conditions for the rest of Taylor coefficients. At our observation point $\x=(0,y^*)$, we have that 
$\partial_x^j = \partial_n^j$, since the $n$-axis is parallel to the $x$-axis. 
We can therefore easily transform the solutions in the ray-centered coordinates $t,n$ 
back to the coordinate system $x, s$. In the end we note that all terms
with odd $x$-derivatives are zero. Hence,
we obtain that $a_1(0) = b_1(0)=0$ and 
$e_{12}+e_{22}$ in \eq{e1def} and \eq{e2def} simplifies to
$$
e_{12}(0)+e_{22}(0) = m'(0) \, a_2(0) \, d_2(0) + m'(0) \, A(0,0) \, b_2(0) \, d_4(0).
$$
After some additional
algebraic manipulations and assuming that $P_0 = i$, $\Im Q_0=0$, $\Re Q_0 >0$, we get
\begin{align}
\label{a2_0}
a_2(0) &= i \frac{3 Q_0^{1/2} \, y^* - 2 Q_0^{1/2} \, (Q_0+i y^*)^2 m'(0) \frac{d}{ds} \theta(y^*)}{4(Q_0+i y^*)^{7/2}}, \\
\label{b2_0}
b_2(0) &= i\frac{(Q_0^2 + {y^*}^2) \left( -y^* +4 (Q_0+i y^*)^2 m'(0) \frac{d}{ds} \theta(y^*)\right)}{8 Q_0 (Q_0+i y^*)^4 },\\
\label{c1_0}
c_1(0) &= \frac{y^* + (Q_0^2 + {y^*}^2) m'(0) \frac{d}{ds} \theta(y^*)}{2 Q_0}, 
\end{align}
and
\begin{equation}
A(0,0) = \frac{Q_0^{1/2}}{(Q_0+i y^*)^{1/2}},
\qquad
\wi(0) = \left( \frac{Q_0 ^2 + {y^*}^2}{\omega \, Q_0} \right)^{1/2}.
\end{equation}
Moreover, by (\ref{tet_const}-\ref{t_const}),
\begin{equation}\label{m_prime}
\frac{d}{ds} \theta(y^*) = y_0''(0), \quad m'(0)={(X^{-1})}'(0) = \left( 1- y^* y_0''(0)\right)^{-1}.
\end{equation}
Therefore, knowing $y_0(s)$ and by (\ref{a2_0}-\ref{m_prime}) and (\ref{d_p}), we can calculate 
$$e_{12}(0)+e_{22}(0)=e^{-i\omega\phi(0,m(0))}C^*(0)=e^{-i\omega y^*}C^*(0). 
$$
Note that $C^*(0)$ only depends on $Q_0$, $y^*$ and $y_0''(0)$. 

We now consider the following two canonical cases:
\begin{itemize}
\item [(1)]  $y_0''(0) = 0$, 
\item [(2)]  $y_0''(0) = -1$.
\end{itemize}
The first case corresponds to a line $y_0=0$. The second case corresponds to a circle $y_0(s)= -1+\sqrt{1 -s^2}$ or a parabola $y_0(s)= -s^2/2$. Note that with an initial curve with positive second derivative, the rays will intersect and form a caustic, and then our theory does not hold.

For the first case, we obtain the simple expression
\be{E_even_1}
C^*_{\rm l}(0) = e^{-i\omega y^*}\frac{n_0 \, y^* \, Q_0^2}{(Q_0 + i y^*)^2 \, (Q_0^2 + {y^*}^2)^{3/2}},
\ee
where $n_0$ is a constant complex number.
For the second case, the expression is much more complicated. 
In the small and large $Q_0$-limit we have
\be{E_even_2}
C^*_{\rm c}(0) =e^{-i\omega y^*}\frac{n_1+n_2y^*+n_3{y^*}^2}
{{y^*}^{4}\sqrt{1+y^*}}Q_0^{2}+\mathcal{O}(Q_0^3),\qquad
C^*_{\rm c}(0) =e^{-i\omega y^*}\frac{n_4}{\sqrt{1+y^*}}Q_0^{-5/2}+\mathcal{O}(Q_0^{-7/2}),
\ee
where $n_j$, with $j=1, \dotsc, 4$, are constant complex numbers. 
The amplitude of the geometrical optics solution is proportional to $|1-y^* \, y_0''(0)|^{-1/2}$, and by \eq{nonoscerror}, the 
relative error will be
$$
|E_{\rm rel}| = |E^T| \, |1-y^* \, y_0''(0)|^{1/2} 
+\mathcal{O}(\omega^{-3/2})
= \omega^{1/2} {\wi}^3(0) \, |C^*(0)| \,|1-y^* \, y_0''(0)|^{1/2}
+\mathcal{O}(\omega^{-3/2}).
$$
We therefore obtain the leading order term
$$
|E_{\rm rel}^{\rm l}| =\omega^{-1}\;\left| \frac{n_0 \, y^* Q_0^{1/2}}{(Q_0 + i y^*)^2} \right|,
$$
and for small and large $Q_0$,
$$
|E_{\rm rel}^{\rm c}| =\omega^{-1}\,
\frac{n_1+n_2y^*+n_3{y^*}^2}
{{y^*}}\sqrt{Q_0}
+\mathcal{O}(\omega^{-1}Q_0^{3/2}),\qquad
|E_{\rm rel}^{\rm c}| =\omega^{-1}\,
\frac{n_4}{Q_0}+\mathcal{O}(\omega^{-1} Q_0^{-2}),
$$
corresponding to \eq{E_even_1} and \eq{E_even_2}, respectively. 

\begin{figure}[!t]
  \begin{center}
   \subfigure[]{\includegraphics[width=.45\textwidth]{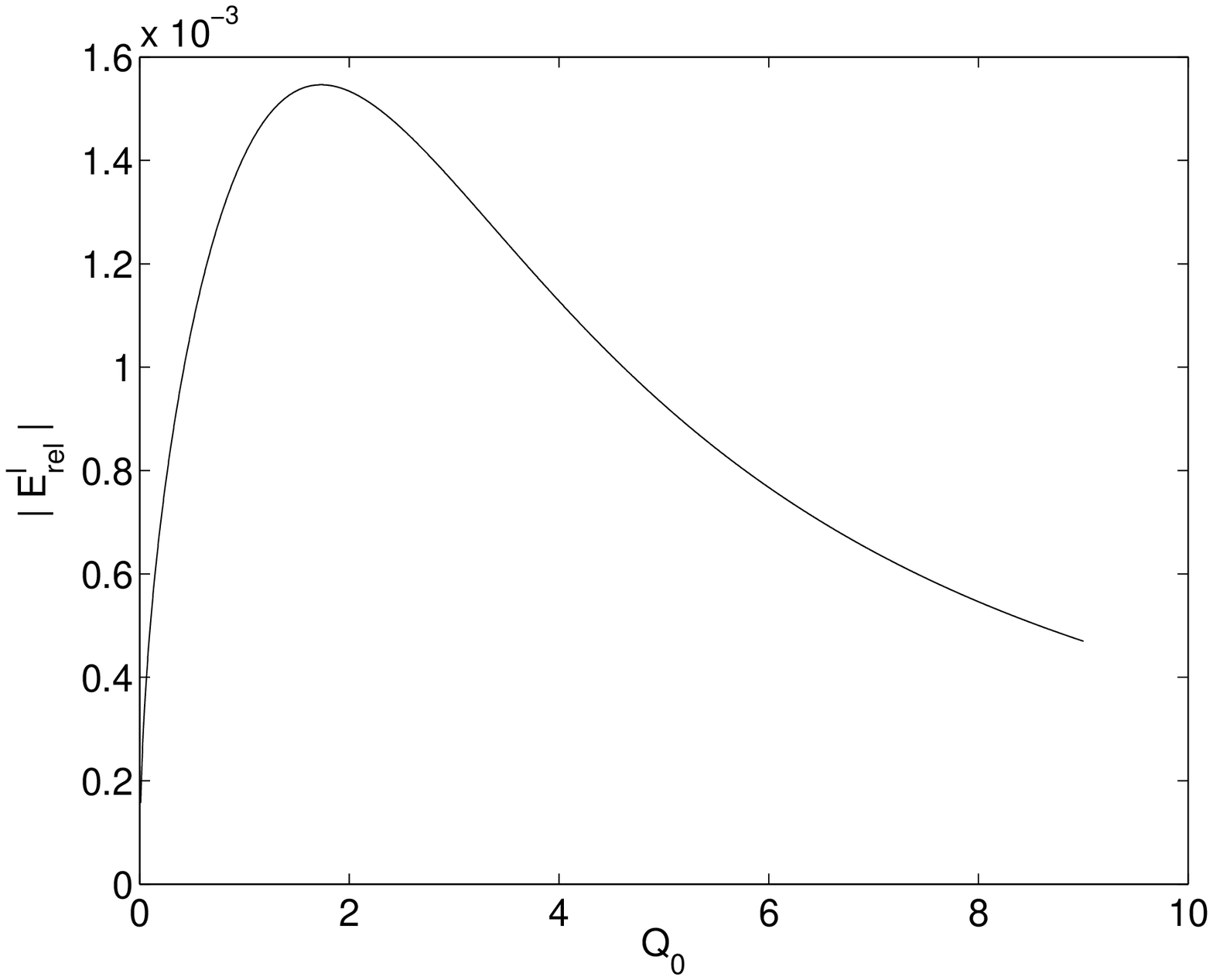}}
   \subfigure[]{\includegraphics[width=.45\textwidth]{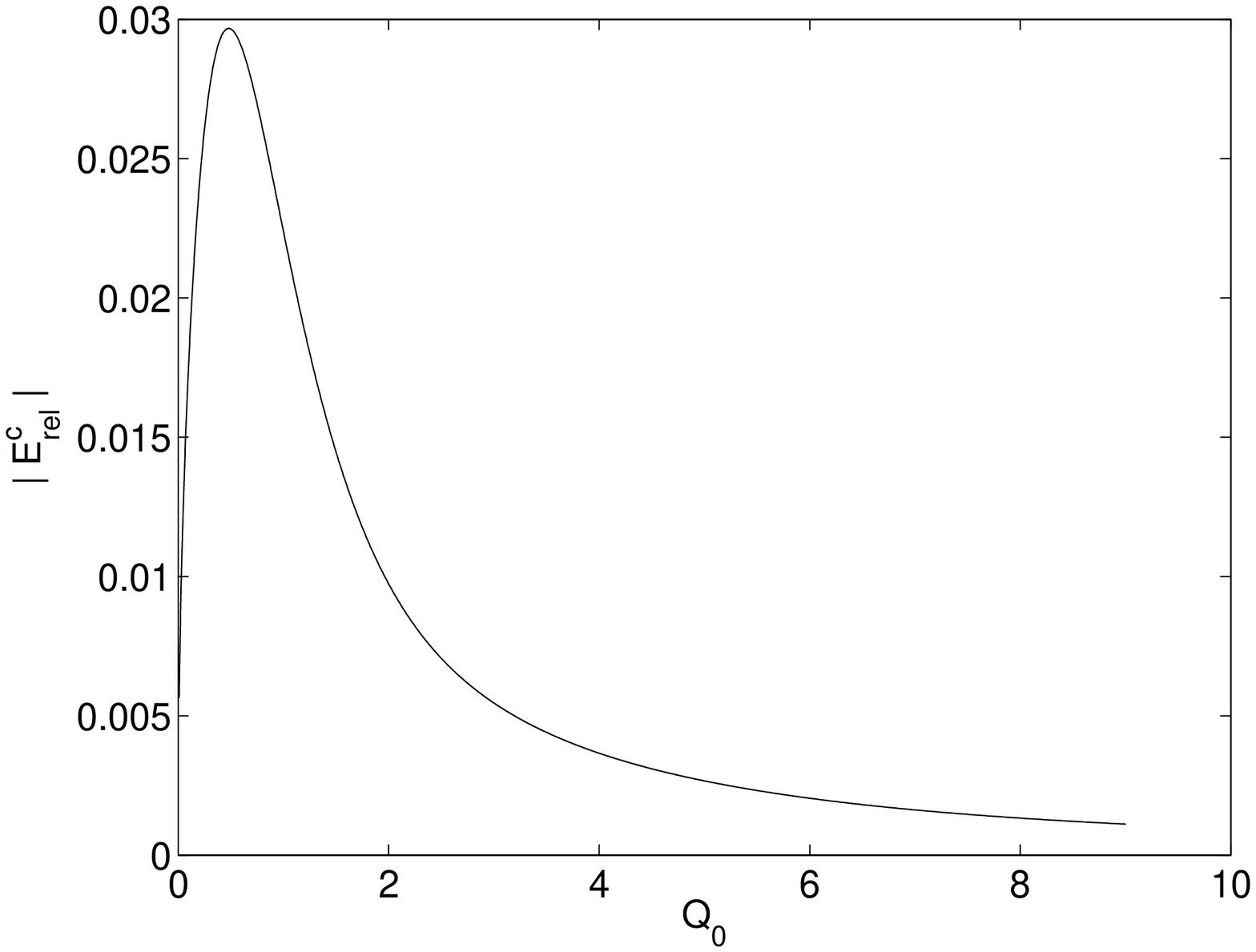}}\\
   \subfigure[]{\includegraphics[width=.45\textwidth]{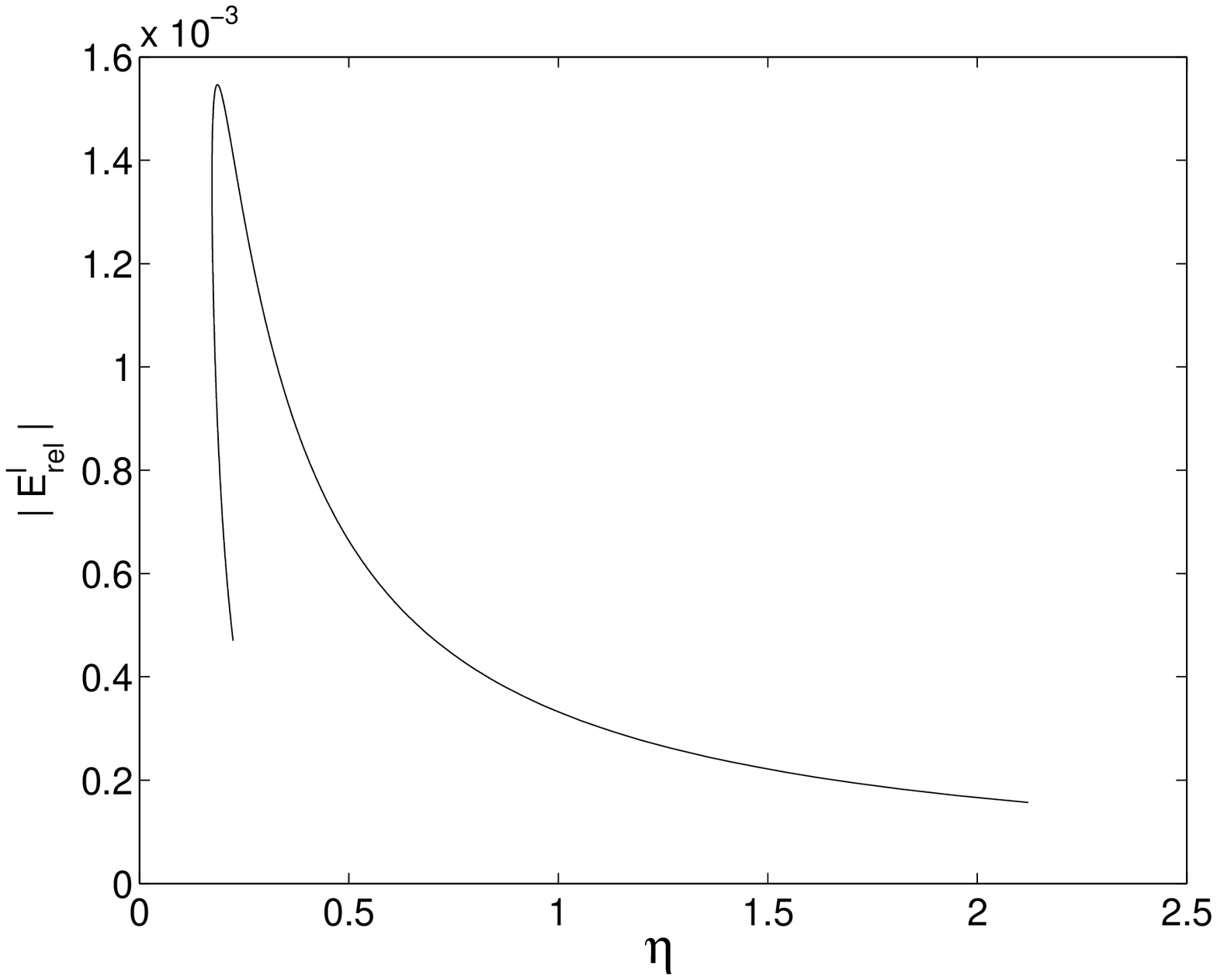}}
   \subfigure[]{\includegraphics[width=.45\textwidth]{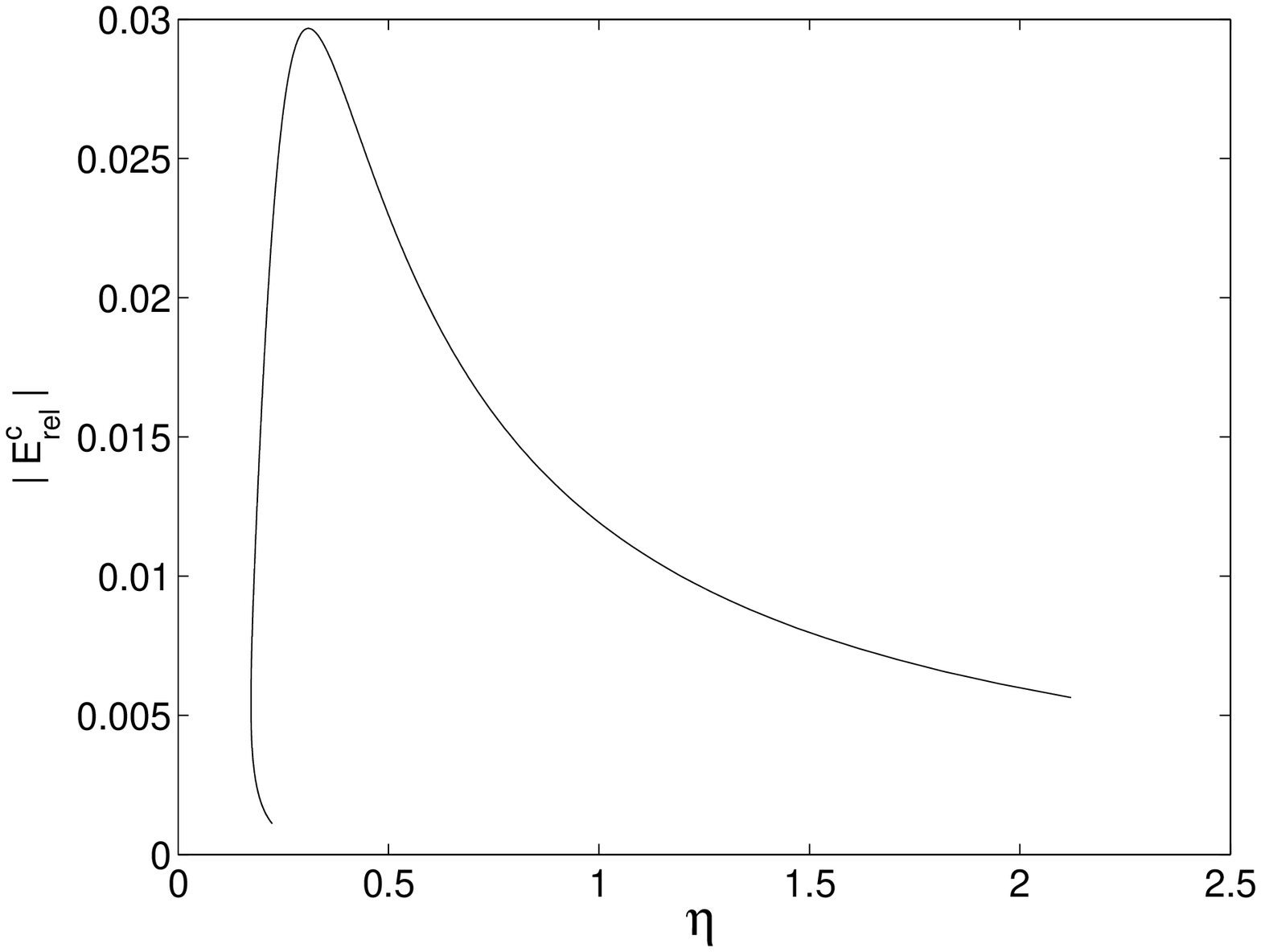}}\\
  \caption{Absolute value of relative error as a function of $Q_0$ (top) and of the width $\wi$ (bottom) in the case when the initial source is a line (left) and a circle (right).}
\lbfig{Rel_Error}
  \end{center}
\end{figure}
\fig{Rel_Error} shows the absolute values of the relative errors at $y^*=3$. 
As it can be seen from the formulas and figures, the relative error 
has a direct relation with $Q_0$, tending to zero both for small and large $Q_0$.
A reduced error with large $Q_0$ has also been noticed in \cite{Klimesh86} for the oscillatory part of the error (or the discretization error). 
However, there is no clear connection between the error and the beam width;
the same width can correspond to different errors.

In many approximations, the optimal $Q_0$, corresponding to the minimum beam width at a receiver point is chosen for computations, see \cite{Cerveny_etal} for instance.
\fig{Q0_W} shows the beam width as a function of $Q_0$. 
In our case the minimum width is attained at $Q_0 = y^*$.
With $Q_0 = y^*$ and $y*\gg 1$ we get
$$ 
|E_{\rm rel}^{\rm l}| = \frac{N}{\omega \, {y^*}^{1/2}}, \qquad |E_{\rm rel}^{\rm c}| \approx \frac{N'}{\omega y^*},
$$
with $N$ and $N'$ being constant numbers. 
When using this $Q_0$, we do not obtain the minimum error as was seen above. However, importantly, the relative error decreases as the distance from the source increases.

We conclude that in the case analysed here large and very small $Q_0$ will improve
the Taylor expansion error.  From \fig{Q0_W} we see that this
corresponds to having {\em wide} beams, not narrow beams. 
One should keep in mind, however,
that this is not the whole story. 
The approximation of the initial data where
the source curve is not flat
and/or the amplitude is not constant
will in general deteriorate
when wider beams are used. Hence, this restricts the beam
widths that can be used. 
Wider beams also mean that the wave field will be more
expensive to evaluate since beams contribute more globally to the solution.
Moreover, our result is strictly for constant coefficients.
In the presence of a varying speed of propagation where the properties may change dramatically as we get farther from the central rays, the Taylor expansion error 
could be large for wide beams. In addition, when the rays can bend, it may not be possible to have very wide beams, since as was noted before, the Gaussian beam approximation may break down when the phase becomes non-smooth, and this happens at some distance away from the central ray (outside {\it the regularity region}).
In the general case, finding the optimal $Q_0$ 
for a given observation point is an open problem.

\begin{figure}[!t]
  \begin{center}
{\includegraphics[width=.6\textwidth]{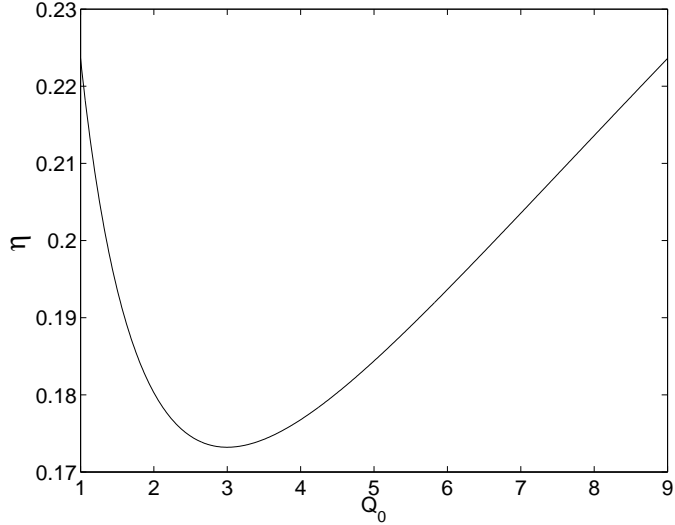}}
\caption{The beam width as a function of $Q_0$ at $y^*=3$.}
\lbfig{Q0_W}
  \end{center}
\end{figure}







\end{document}